\documentclass[a4paper,12pt]{article}
\usepackage{amsmath}
\usepackage{amssymb}
\usepackage{amsthm}
\usepackage{mathtools}
\usepackage{stmaryrd}
\usepackage{bbold}
\usepackage{breqn}
\usepackage{txfonts}
\usepackage{graphicx}
\usepackage{xfrac}

\usepackage[dvipsnames]{xcolor}
\usepackage{soul}
\usepackage{comment}
\specialcomment{coloredcomment}{\begingroup\color{blue!50}}{\endgroup}

\usepackage{accents}

\usepackage[english]{babel}
\usepackage{blindtext}

\usepackage{bm}
\bmdefine{\ba}{a}
\bmdefine{\bb}{b}
\bmdefine{\bc}{c}
\bmdefine{\bd}{d}
\bmdefine{\be}{e}
\bmdefine{\bf}{f}
\bmdefine{\bg}{g}
\bmdefine{\bj}{j}
\bmdefine{\bk}{k}
\bmdefine{\bn}{n}
\bmdefine{\bp}{p}
\bmdefine{\bq}{q}
\bmdefine{\br}{r}
\bmdefine{\bx}{x}
\bmdefine{\by}{y}
\bmdefine{\bu}{u}
\bmdefine{\bv}{v}
\bmdefine{\bw}{w}
\bmdefine{\bz}{z}
\bmdefine{\bE}{E}
\bmdefine{\bB}{B}
\bmdefine{\bC}{C}
\bmdefine{\bD}{D}
\bmdefine{\bJ}{J}
\bmdefine{\bR}{R}
\bmdefine{\bS}{S}
\bmdefine{\bT}{T}
\bmdefine{\bvarphi}{\varphi}
\bmdefine{\bphi}{\phi}
\bmdefine{\bpsi}{\psi}
\bmdefine{\bchi}{\chi}
\bmdefine{\bPsi}{\Psi}
\bmdefine{\blambda}{\lambda}
\bmdefine{\bmu}{\mu}
\bmdefine{\bnu}{\nu}

\def\xfoo#1^#2\relax#3\valign{%
\mathbf{#1}\ifx\valign#2\valign\else^{\mathbf{#2}}\fi}

\usepackage[mathscr]{euscript}

\makeatletter
\newcommand*{\defeq}{\mathrel{\rlap{%
					 \raisebox{0.3ex}{$\m@th\cdot$}}%
					 \raisebox{-0.3ex}{$\m@th\cdot$}}%
					 =}
\makeatother

\makeatletter
\newcommand*{\eqdef}{=\mathrel{\rlap{%
					 \raisebox{0.3ex}{$\m@th\cdot$}}%
					 \raisebox{-0.3ex}{$\m@th\cdot$}}%
					 }
\makeatother

\newcommand{\sgn}{\mathop{\mathrm{sgn}}}

\newcommand{\op}[1]{\operatorname{#1}}
\newcommand{\uo}{\operatorname{\mathfrak{u}}}
\newcommand{\bbb}[1]{\mathbb{#1}}
\newcommand{\llb}{\llbracket}
\newcommand{\rrb}{\rrbracket}

\def\XXint#1#2#3{{\setbox0=\hbox{$#1{#2#3}{\int}$}
     \vcenter{\hbox{$#2#3$}}\kern-.5\wd0}}

\newcommand\wrapped[1]%
  {\renewcommand\arraystretch{1}%
   \begin{array}{@{}l@{}}#1\end{array}%
  }
  
\theoremstyle{definition}
\newtheorem{example}{Example}

\theoremstyle{definition}

\theoremstyle{definition}

\usepackage{mathtools}

\DeclarePairedDelimiterX\braket[2]{\langle}{\rangle}{#1\,\delimsize\vert\,\mathopen{}#2}

\usepackage[]{hyperref}
\hypersetup{
    pdftitle={Nonlinear electromagnetic dispersion relation for a magnetized plasma},
    pdfauthor={Roberto Ricci},
    pdfsubject={Magnetized plasma kinetic model},
    pdfkeywords={magnetized plasma, nonlinear dispersion relation},
    bookmarks=true,
    bookmarksnumbered=true,     
    bookmarksopen=true,         
    colorlinks=true,            
    pdfstartview=Fit,           
    pdfpagemode=UseOutlines,    
    pdfpagelayout=TwoPageRight
}

\usepackage{cleveref}
\crefformat{pluralequation}{#2eqs.~(#1)#3}
\Crefformat{pluralequation}{#2Eqs.~(#1)#3}

\usepackage{authblk}

\begin{document}
	
\title{Umbral theory and the algebra of formal power series}
\author{Roberto Ricci\footnote{E-mail: roberto.ricci@enea.it (Roberto Ricci)}}
\affil[1]{ENEA, Nuclear Department NUC-DTT, Frascati Research Center, Via E. Fermi 45, 00044 Frascati (Rome), Italy}
\date{}
\setcounter{Maxaffil}{0}
\renewcommand\Affilfont{\itshape\small}
\maketitle

\section*{Abstract}

Umbral theory, formulated in its modern version by S. Roman and G.~C. Rota, has been reconsidered in more recent times by G. Dattoli and collaborators with the aim of devising a working computational tool in the framework of special function theory. Concepts like umbral image and umbral vacuum have been introduced as pivotal elements of the discussion, which, albeit effective, lacks of generality.

This article is directed towards endowing the formalism with a rigorous formulation within the context of the formal power series with complex coefficients $(\bbb{C}\llb t \rrb, \partial)$. The new formulation is founded on the definition of the umbral operator $\uo$ as a functional in the "umbral ground state" subalgebra of analytically convergent formal series $\varphi \in \bbb{C}\{t\}$. 

We consider in detail some specific classes of umbral ground states $\varphi$ and analyse the conditions for analytic convergence of the corresponding umbral identities, defined as formal series resulting from the action on $\varphi$ of operators of the form $f(\zeta \uo^\mu)$ with $f \in \bbb{C}\{t\}$ and $\mu, \zeta \in \bbb{C}$. For these umbral states, we exploit the Gevrey classification of formal power series to establish a connection with the theory of Borel-Laplace resummation, enabling to make rigorous sense of a large class of -- even divergent -- umbral identities. 

As an application of the proposed theoretical framework, we introduce and investigate the properties of new umbral images for the Gaussian trigonometric functions, which emphasise the trigonometric-like nature of these functions and enable to define the concept of "Gaussian Fourier transform", a potentially powerful tool for applications.

\section*{Introduction}\label{sec:introduction}

Classical umbral calculus is a symbolic technique for manipulating polynomial sequences, along with associated  generating functions, in a way that resembles formal substitution. 

The term\footnote{It was coined by J.~J. Sylvester from the Latin word "umbra", meaning "shadow".} reflects the fact that, for many types of identities involving sequences of polynomials with powers $a^n$, "shadow" identities are obtained when the polynomials are changed to discrete values and $a^n$ is changed to the falling Pochhammer symbol $(a)_n=a(a-1)...(a-n+1)$ \cite{aDB00}.

\noindent "Umbral calculus" is currently mainly used to refer to the study of Sheffer sequences\footnote{S. Roman himself describes umbral calculus as the study of the class of Sheffer sequences \cite[p. 2]{sR84}.}, including polynomial sequences of binomial type and Appell sequences, often by means of techniques borrowed from the calculus of finite differences. Its origins trace back to the 19th century \cite{jB1861}, but the first rigorous framework was introduced in the 1970's by Gian-Carlo Rota and collaborators \cite{gcR75}, and later systematized and extended by Steven Roman \cite{sR84}.

In his seminal work, Rota sought to place the classical umbral methods of symbolic computation on solid algebraic foundations. The key insight in Rota's theory is to reinterpret the manipulation of sequences and polynomials in terms of linear functionals acting on the ring of polynomials $\mathbb{C}[x]$, thus connecting umbral calculus with combinatorics.

Roman extended Rota's ideas in a more symbolic and axiomatic direction, culminating in a comprehensive treatment in \emph{The Umbral Calculus}~\cite{sR84}. In particular, Roman formalised the notion of \emph{umbral operators} $\uo$ acting on formal power series and introduced the concept of \emph{ground states} $\varphi$ -- special formal series representing the identity of evaluation.

While Rota's formulation is algebraically motivated and deeply connected to combinatorics, Roman's version of umbral calculus is more formal and symbolic, enabling the construction of generalised umbral transformations, especially in the analysis of special functions.
Both frameworks have found broad applications in:
\begin{itemize}
    \item combinatorics and enumerative identities;
    \item operational calculus and generating functions;
    \item asymptotic expansions and formal series;
    \item special function theory (e.g. hypergeometric functions, orthogonal polynomials).
\end{itemize}

A different declination of umbral calculus has been recently proposed by G. Dattoli and collaborators (see \cite{sL22} for a thorough introduction and \cite{gD24} for some recent results) with the specific aim of formulating a working computational tool in the framework of special function theory. In their investigations, concepts like "umbral image" and "umbral vacuum" are used as pivotal elements, although the corresponding mathematical entities are usually introduced at an informal level, lacking a rigorous definition.

This version of umbral calculus, sometimes referred to as "indicial umbral calculus", is devised as an effective method to simplify complex calculations through the identification of "umbral identities" among some highly transcendental functions and their \emph{umbrae}, i.e. some simpler operator functions acting on "ground states" defined \emph{ad hoc}. A typical umbral identity is an expression of the form 
\begin{dmath}\label{eq:umbral_id}
	F(z) = f(z\uo^{\,\mu})\, \bvarphi
\end{dmath},
where $F$ and $f$ are functions expandable in Maclaurin series and $\mu$ is a complex parameter, interpreted as a power exponent of the "umbral operator" $\uo$, eventually acting on the ground state $\bvarphi$.

Once a relation like \cref{eq:umbral_id} has been established, any calculation involving $F(z)$ is carried out by using $f(z\uo^{\,\mu})$ in its place and treating $\uo$ as an ordinary numeric parameter. The outcome of the calculation is eventually applied to $\bvarphi$, yielding the final result in the form of a (usually converging) series expansion. The validity of the procedure is verified \emph{ex-post}, by confronting the results obtained in this way with those attainable with more traditional, analytic or numerical methods.  

\noindent Although lacking a general justification, the afore described method is highly effective and mostly provides correct results with minimum effort. Nonetheless, in some cases it fails to produce convergent expansions, with no clear explanation.

In the present work, a proposal is made for a new general approach aimed at clarifying: 1) why and how indicial umbral calculus works; 2) in which conditions it apparently fails; 3) what can be done, in this case, in order to recover meaningful results.

\noindent The proposed approach provides a solid foundation within the differential algebra of formal power series with complex coefficients $(\bbb{C}\llb t \rrb, \partial)$ \cite{dS16}. Specifically, the umbral operator $\uo$ is rigorously defined as a functional in the "umbral ground state" subalgebra of analytically convergent formal series $\bbb{C}\{t\} \subset \bbb{C}\llb t\rrb$. 

\noindent After singling out some specific classes of umbral ground states $\varphi \in \bbb{C}\{t\}$, which comprise the vast majority of the cases considered so far in the literature on the topic, we investigate umbral identities of the form $F(\zeta) = f(\zeta \uo^\mu)[\varphi]$. These are defined as formal series in the isomorphic differential algebra $\bbb{C}\llb \zeta\rrb$, resulting from the action on $\varphi$ of operators of the form $f(\zeta \uo^\mu)$ with $f \in \bbb{C}\{t\}$ and $\mu, \zeta \in \bbb{C}$. 

With an eye to applications, where one is interested to the case where $F(\zeta)$ is a full-fledged function of the complex variable $\zeta$, we analyse the conditions for analytic convergence of these formal series. To this aim, we make recourse to the well-established theory of Gevrey classification of formal series in nested subalgebras of increasing level -- that is, roughly speaking, increasing "divergence order". Once characterised a certain divergent formal series $F$ as belonging to a certain Gevrey subalgebra of index $k>0$, we make use of the well-established theory of k-Borel-Laplace resummation in order to find the $k$-sum of $F$, i.e. a function that is $k$-Gevrey asymptotic to $F$ in an open sector of the complex plain.

The idea of using the Borel transform as a means to provide a solid foundation for the umbral method is not new. In particular, the role of integral transforms of the Borel type as a bridge between the heuristic umbral approach and analytical techniques was investigated by G. Dattoli and collaborators in a series of articles \cite{gD15, gDsL20, gD20}. The present proposal differentiates from those early attempts in suggesting a systematic and theoretically well-founded recourse to the \emph{formal} Borel transform in the wider context of formal series.

\vskip 0.3cm

The article is organised as follows.

\noindent In \cref{sec:math_framework} we provide a brief introduction to the algebra of formal power series in one indeterminate $\bbb{C}\llb t \rrb$. In particular we informally review the main results of the theories of Gevrey classification and asymptotics and Borel-Laplace resummation that are used in the rest of the paper. 

\noindent In \cref{sec:umbral_theory} we present our proposal for a rigorously founded indicial umbral theory in the context of formal power series, by providing sound definitions of the umbral operator $\uo$ and of umbral identities, reinterpreted as formal series resulting from the action on $\varphi$ of well-defined operators of the form $f(\zeta \uo^\mu)$.

\noindent In \cref{sec:umbral_ground_states} we consider some specific classes of umbral ground states and analyse the analytic convergence conditions of umbral identities relying on such ground states, by exploiting the mathematical framework outlined in \cref{sec:formal_series_algebra}.

\noindent In \cref{sec:gaussian_trig} we suggest a new umbral formulation of the Gaussian trigonometric functions as a case study for the proposed framework. The new formulation  emphasises the trigonometric-like nature of these functions and enables to define the concept of "Gaussian Fourier transform".

\noindent In \cref{sec:conclusions} we draw some conclusions -- including a brief critical comparison of the framework proposed in this article with Roman's classical umbral calculus -- and expose ideas for future developments. 

We use the following standard notations: $\bbb{N}$ is the set of natural numbers; $\bbb{N}^0 = \{0\} \cup \bbb{N}$, the set of non-negative integer numbers; $\bbb{Z}^-$ is the complementary set of negative integer numbers; $\bbb{R}^+$ is the set of strictly positive real numbers.

\section{The mathematical framework}\label{sec:math_framework}

We give here a brief, mostly informal exposition of the mathematical framework used in the following sections.

\subsection{The algebra  $\mathbb{C}\llb t\rrb$}\label{sec:formal_series_algebra}

Let us consider the set of all formal power series with complex coefficients in the indeterminate $t$ \cite{aD15,dS16}:
\begin{dmath*}
	\bbb{C}\llb t \rrb \defeq {
		\left\{\varphi(t) = \sum_{n = 0}^\infty c^{(\varphi)}_n t^n \,\middle|\, c^{(\varphi)}_n \in \bbb{C}\right\}
	}
\end{dmath*}.
Formal series can be summed and multiplied by complex numbers term by term, so they form a complex vector space. 

Given two formal series $\varphi$ and $\chi$, their product $\xi = \varphi\chi$ realises another formal series with coefficients:
\begin{dmath*}
	c^{(\xi)}_n \defeq \sum_{k = 0}^n c^{(\varphi)}_k c^{(\chi)}_{n-k} 
\end{dmath*},
which are the result of a Cauchy discrete convolution \cite{cC15}.

\noindent This makes of $\bbb{C}\llb t \rrb$ an algebra.
The corresponding differential algebra $(\bbb{C}\llb t \rrb, \partial)$ is obtained after equipping $\bbb{C}\llb t \rrb$ with he additional structure provided by $\partial \defeq \frac{\op{d}}{\op{d}t}$ -- the natural derivation with respect to the indeterminate, satisfying the usual Leibniz rule.

We define the valuation of $\varphi$ as the value of $n$ corresponding to the first non-zero coefficient, so $\op{val}(0) = \infty$, $\op{val}(\varphi) = 0$ if $\varphi$ has nonzero constant term  and, in general, $\op{val}(\varphi) \geq k$ if $\varphi \in t^k\bbb{C}\llb t \rrb$, where:
\begin{dmath*}
	t^k\bbb{C}\llb t \rrb \defeq {
		\left\{\varphi(t) = \sum_{n = k}^\infty c^{(\varphi)}_n t^n \,\middle|\, c^{(\varphi)}_n \in \bbb{C}\right\} \condition*{k \in \bbb{N}^0}
	}
\end{dmath*}.
Using the valuation it is possible to define the quantity:
\begin{dmath*}
	\op{d}(\varphi, \chi) \defeq 2^{-\op{val}(\varphi - \chi)}
\end{dmath*},
which can be easily recognised as a distance function, since it meets the four axioms \cite{dB01}:
\begin{enumerate}
	\item non-negativity: $\op{d}(\varphi, \chi)\ge 0$;
	\item identity of indiscernibles: $\op{d}(\varphi, \chi) =0 \iff \varphi=\chi$;
	\item symmetry: $\op{d}(\varphi, \chi) = \op{d}(\chi, \varphi)$;
	\item triangle inequality: $\op{d}(\varphi, \chi)\le \op{d}(\varphi,\xi)+\op{d}(\xi,\chi)$,
\end{enumerate}
for any $\varphi, \chi, \xi \in \bbb{C}\llb t \rrb$. 

The topology induced by $\op{d}$, known as Krull topology \cite{dS16}, makes of the algebra $\bbb{C}\llb t \rrb$ a complete metric space and enables to introduce the concept of formal convergence of formal series. In particular, a sequence of formal series $\{\varphi_p\}_{p=0}^\infty$ is a Cauchy sequence and converges formally to $\varphi$ if and only if, for each $n \in \bbb{N}^0$, the sequence of n-th coefficients $\{c^{(\varphi_p)}_n\}_{p=0}^\infty$ is stationary, that is for each $n$ there is an integer $\mu(n)$ such that the subsequence $\{c^{(\varphi_p)}_n\}_{p=\mu(n)}^\infty$ is constant with $c^{(\varphi_p)}_n = c^{(\varphi)}_n$, for $p\geq \mu(n)$. Then, obviously:
\begin{dmath*}
	\lim_{p \to \infty} \varphi_p(t) = {\varphi(t) = \sum_{n=0}^\infty c^{(\varphi)}_n t^n}
\end{dmath*}.
Since a series of formal series is defined as the sequence of its partial sums, which are themselves formal series,
\begin{dmath*}
	\sum_{p=0}^\infty \varphi_p(t) = {\{\sigma_n(t)\}_{n=0}^\infty, \;\;\; \sigma_n(t) \defeq \sum_{p=0}^n \varphi_p(t)}
\end{dmath*},
we have the following practical criterion of formal convergence for a series of formal series: each term $\varphi_p$ must belong to $t^{\nu_p}\bbb{C}\llb t \rrb$, that is $c^{(\varphi_p)}_n = 0$ for $n < \nu_p$, for some integer $\nu_p$ such that $\lim_{p \to \infty} \nu_p = \infty$. Then the sum converges to the formal series $\varphi$ with coefficients
\begin{dmath*}
	c^{(\varphi)}_n =
		\sum_{p \in M_n} c^{(\varphi)_p}_n \condition*{M_n = \{p \,|\, \nu_p \leq n\}}
\end{dmath*}.

It is important to stress that formal convergence and analytic convergence are totally unrelated concepts. Formal convergence has to do with rules ensuring that the result of certain manipulations of formal series is itself a well-formed formal series. In this sense, any formal series is by definition formally convergent. Needless to say, this is not true for standard analytic convergence of a generic formal series when interpreted as an ordinary power series in $\bbb{C}$.

\subsection{Gevrey classification and asymptotics}

\subsubsection{Gevrey classification of elements of $\mathbb{C}\llb t \rrb$}

 The Gevrey classification of the elements of $\mathbb{C}\llb t \rrb$ \cite{aD15,bP16} is of fundamental importance for the considerations that will be made in the following.  
 
 A formal series $\varphi$ is said to be Gevrey of level $k > 0$ ($k$-Gevrey for short) if some constants $C, A >0$ exist such that
\begin{dmath}\label{eq:k-Gevrey_condition}
		| c^{(\varphi)}_r | \leq C A^r (r!)^{\frac{1}{k}}\condition{for any $r\in \bbb{N}^0$}
\end{dmath}.
It follows from the definition that a formal series $\sum_{r=0}^\infty c^{(\varphi)}_r t^r$ is $k$-Gevrey if and only if the series $\tilde{\varphi} = \sum_{r=0}^\infty c^{(\varphi)}_r/(r!)^{\frac{1}{k}} t^r$ is convergent in the analytic sense.

For each $k$, we denote  with $\mathbb{C}\llb t \rrb_{\frac{1}{k}}$ the subalgebra formed by all $k$-Gevrey formal series.

\noindent These subalgebras are nested, in the sense that
\begin{dmath}\label{eq:gevrey_nested}
	\bbb{C}\{t\} = {
		\bbb{C}\llb t \rrb_0 \subset 
		\bbb{C}\llb t \rrb_{\frac{1}{k}} \subset 
		\bbb{C}\llb t \rrb_{\frac{1}{k'}} \subset
		\bbb{C}\llb t \rrb_\infty =
		\bbb{C}\llb t \rrb
	}
\end{dmath},
for any $k, k'$ satisfying $\infty > k > k' > 0$.

\noindent So, a $k$-Gevrey series is also $k'$-Gevrey for any $k'<k$. A formal series is \emph{exactly} $k$-Gevrey if $k$ is the greatest value for which \cref{eq:k-Gevrey_condition} is satisfied. If such value of $k$ is infinite, i.e. $| c^{(\varphi)}_r | \leq C A^r$, the corresponding subalgebra $ \mathbb{C}\llb t \rrb_0$ coincides with the subalgebra of analytically convergent series, which is denoted with $\bbb{C}\{t\}$. 
So, a 0-Gevrey series $\varphi$ converges analytically in an open disc $D(R_\varphi)$ centred in the origin, with (possibly infinite) radius $R_\varphi$.

\subsubsection{Gevrey asymptotics}

We denote with $S(d, \Theta, \rho)$ an open sector of the complex plain (or of the Riemann surface of the complex logarithm if $\Theta > 2\pi$) centred in the origin, bisected by the half-line $d$, of aperture $\Theta$ and radius $\rho > 0$. We use the simplified notations  $S(d, \Theta)$ for the infinite sector with infinite radius, and  $S(\Theta, \rho)$ for the sector bisected by the positive real axis $\bbb{R}^+$.
 
Let $V = S(d, \Theta, \rho)$ be an open sector as before, $s$ a function holomorphic in $V$, $\varphi \in \mathbb{C}\llb t \rrb$ a formal series and $k>0$. We say that $\varphi$ is the $k$-Gevrey asymptotic ($k$-asymptotic for short) expansion of $s$ in $V$ and write
\begin{dmath*}
	\varphi \sim_{\overset{1/k}{V}} s 
\end{dmath*}
if, for any $n \in \bbb{N}^0$ and any compact subsector $W$ of $\{0\} \cup V$, it is possible to find two constants $A, C >0$ such that
\begin{dmath*}
	\Bigg\lvert s(t) - \sum_{r=0}^{n-1} c^{(\varphi)}\,t^n \Bigg\rvert \leq C A^n (n!)^{\frac{1}{k}} |\,t^n\,| 
\end{dmath*}.

Note that if $\varphi$ is the $k$-Gevrey asymptotic expansion of some function $s$ holomorphic in $V$, then $\varphi \in \mathbb{C}\llb t \rrb_{\frac{1}{k}}$, i.e. it is $k$-Gevrey. As a consequence, if $A_k(V)$ is the differential algebra comprising all the functions holomorphic in $V$ that admit there a $k$-Gevrey asymptotic expansion, it is possible to define a differential algebra homomorphism $\mathcal{I}_k:A_k(v) \to \mathbb{C}\llb t \rrb_{\frac{1}{k}}$, which, at each function $s \in A_k(V)$, associates a $\varphi \in \mathbb{C}\llb t \rrb_{\frac{1}{k}}$ such that $\varphi$ is a $k$-asymptotic expansion of $s$ in $V$. Namely:
\begin{dmath*}
	s \mapsto {\mathcal{I}_k(s) = \varphi}\condition*{\varphi \sim_{\overset{1/k}{V}} s}
\end{dmath*}.
It is natural to wonder whether such homomorphism is injective (different functions have different $k$-asymptotic developments), surjective (any $k$-Gevrey formal series is the $k$-asymptotic development of some function) or both. 

It is possible to demonstrate \cite{aD15} that, if $\varphi \in \mathbb{C}\llb t \rrb_{\frac{1}{k}}$ and $V = S(d, \Theta, \rho)$ is an open sector with $\Theta < \pi/k$, then there is at least a function $s$ holomorphic in $V$ such that $\varphi \sim_{\overset{1/k}{V}} s$. In other terms, the homomorphism $\mathcal{I}_k$ is surjective if $V$ is a sector of "small" central angle, but it fails to be injective, in the sense that different functions can have the same $k$-asymptotic expansion. In fact, it can be shown that the kernel of $\mathcal{I}_k$ is generated by functions like $\op{e}^{-a/t^{1/k}}$ with $a>0$. This is an undesirable property in most applications, e.g. when $\varphi$ is the formal solution of an ordinary differential equation and we want to determine the single function solution of the same equation.

Conversely, a different theorem \cite{gnW12} proves that, when $\Theta > \pi/k$, if there is a function $s$ such that $\varphi \sim_{\overset{1/k}{V}} s$, then that function is unique. Stated differently, the homomorphism becomes injective in these conditions, but it is no more surjective, so that there are $k$-Gevrey formal series that are not $k$-asymptotic expansions of any function.

A formal series $\varphi$ such that there exists a function $s$ holomorphic in $V = S(d, \Theta, \rho)$ with $\Theta > \pi/k$ verifying the condition $\varphi \sim_{\overset{1/k}{V}} s$ is called $k$-summable and the function $s$ is called 
its $k$-sum in the direction $d$.
The problem of finding the $k$-sum of a $k$-summable formal series by means of the Borel-Laplace resummation procedure is analysed in the following subsection.

\subsection{Borel-Laplace resummation}\label{sec:borel_laplace_resummation}

We preliminary give some definitions, that will be extensively used afterwards.

\subsubsection{Formal Borel transforms}

\noindent Given a formal series $\varphi \in \mathbb{C}\llb t \rrb$, let $\tilde{\varphi} \in t\mathbb{C}\llb t \rrb$ be the  series without constant term obtained from $\varphi$ by a shift of the coefficients:
\begin{dmath*}
	\tilde{\varphi}(t) = {
	t\,\varphi(t) = 
		\sum_{r=0}^\infty c^{(\varphi)}_r\,t^{r+1}
	}
\end{dmath*}.
The formal Borel transform of order $k>0$ ($k$-Borel tranform, for short) of $\tilde{\varphi}$ is the formal series $\tilde{\varphi}_{\op{B}_k}$ in the indeterminate\footnote{
	Although not strictly necessary, it is advisable to consider the Borel transform as an element of an algebra $\mathbb{C}\llb \zeta \rrb$ isomorphic to $\mathbb{C}\llb t \rrb$.
} $\zeta$ defined by \cite{eB88, jpR93}:
\begin{dmath}\label{eq:rigorous_borel_transform_k}
	\tilde{\varphi}_{\op{B}_k}(\zeta) = {
		\op{B}_k [\tilde{\varphi}](\zeta) \defeq \sum_{r=0}^\infty \frac{c^{(\varphi)}_r}{\Gamma(1+r/k)} \,\zeta^r \condition*{k>0}
	}
\end{dmath}.

The transform \labelcref{eq:rigorous_borel_transform_k} can be seen as the result of the action of the linear operator $\op{B}_k : t\mathbb{C}\llb t \rrb \to \mathbb{C}\llb \zeta \rrb$.

\Cref{eq:rigorous_borel_transform_k} can be generalised by replacing the 1 appearing in the argument of the Gamma function with a parameter $\beta \in \bbb{C}$. This change produces the generalised formal Borel transform of order $k$ and shift $\beta$, namely \cite{dS16, mLR16}:
\begin{dmath}\label{eq:rigorous_borel_transform_k_beta}
	\tilde{\varphi}_{\op{B}_{k,\, \beta}}(\zeta) = {
		\op{B}_{k,\, \beta} [\tilde{\varphi}](\zeta) \defeq \sum_{r=0}^\infty \frac{c^{(\varphi)}_r}{\Gamma(\beta + r/k)} \,\zeta^r
	}
\end{dmath},
This is a well defined formal series for any $k>0$ and $\beta \in \bbb{C}$ satisfying the condition\footnote{This condition is justified by considering that, in correspondence with the poles of the Gamma function, i.e. when $r/k + \beta \in \{ 0, -1, -2, \dots \}$, the coefficients of $\tilde{\varphi}_{\op{B}_{k\,\beta}}$ are undefined.
} $r/k + \beta \notin \{0\} \cup \bbb{Z}^-$ for all $r \in \mathbb{N}^0$.

\subsubsection{Laplace transforms}

Before introducing the generalised Laplace transforms entering the theory of Borel-Laplace resummation, we need to fix the class of functions these transforms can be applied to. 

\noindent This class, which we denote with $L(V)$, generally consists of functions $P$ that are holomorphic (at least) in an infinite open sector $V = S(d,\Theta)$ and  integrable in a neighbourhood of the origin. These functions must additionally grow at most exponentially at infinity along $d$, the line bisecting $V$. This means that, for any such $P \in L(V)$ there must exist constants $C,T > 0$ such that the inequality
\begin{dmath}\label{eq:exponential_grow}
	|P(\zeta)|	\leq C\op{e}^{|\zeta|/T}
\end{dmath}
holds for $\zeta \in V$ in the limit $|\zeta| \to \infty$.

Let $P \in L(V)$. Then the Laplace transform of order $k$ along $d$ of $P$ is defined by the integral \cite{aD15, dS16, mLR16}:
\begin{dmath}\label{eq:laplace-transform_k_d}
		\mathcal{L}_k[P](t) = \frac{k}{t^{\,k-1}}\int_d \mathrm{d}\zeta\,P(\zeta)\,\zeta^{k - 1}\, \op{e}^{-(\zeta / t)^k}
\end{dmath}.
The exponential kernel $\op{e}^{-(\zeta/t)^k}$ suppresses the exponential growth of $P(\zeta)$ as long as $|\zeta/t|^k > 1$ and the integral converges absolutely in the so-called "Borel disk" \cite{aD15}, an open disk passing through the origin, bisected by $d$ and of diameter $T$, where $T$ is the positive constant entering \cref{eq:exponential_grow}.

We will be mostly interested to the case where $d$ coincides with the positive real axis $\bbb{R}^+$. In this case \cref{eq:laplace-transform_k_d} obviously becomes the improper integral:
\begin{dmath}\label{eq:laplace-transform_k}
		\mathcal{L}_k[P](t) = \frac{k}{t^{\,k-1}}\int_0^\infty \mathrm{d}\zeta\,P(\zeta)\,\zeta^{k - 1}\, \op{e}^{-(\zeta / t)^k}
\end{dmath}.
Note that $\mathcal{L}_1[P](t)$, under the change of integration variable $\zeta \mapsto s = \zeta/t$, coincides with the ordinary Laplace transform.

A further generalisation of \cref{eq:laplace-transform_k_d} is the Laplace transform of order $k$ and shift $\beta$ along $d$, defined as \cite{aD15}:
\begin{dmath}\label{eq:laplace-transform_k_beta}
	\mathcal{L}_{k,\,\beta}[P](t)
= \frac{k}{t^{\,\beta k-1}} \int_d \mathrm{d}\zeta\,P(\zeta)\,\zeta^{\beta k - 1} \op{e}^{-(\zeta/t)^k}
\end{dmath}, 
which reduces as expected to \cref{eq:laplace-transform_k_d} for $\beta = 1$.

\noindent Despite the presence of the shift $\beta$, the region of convergence of the integral remains the same, because the asymptotic balance at infinity is still governed by the exponential damping from $\op{e}^{-(\zeta/t)^k}$ versus the exponential growth of $P(\zeta)$. So, if $P \in L(V)$, then $\mathcal{L}_{k,\,\beta}[P](t)$ converges in the same Borel disk as $\mathcal{L}_k[P](t)$. The additional powers of $\zeta$ and $1/t$ alter the functional form and analytic structure of the transform, but not its convergence domain. In particular, the prefactor $1/t^{\,\beta k -1}$ does affect the analytic behaviour of the output function, by introducing a branch point at $t = 0$ if $\beta k \notin \mathbb{Z}$.

\subsubsection{$k$-Borel-Laplace sum of a $k$-Gevrey formal series}

The notion of Gevrey level applies to any formal series and is in principle unrelated to the notion of analytic convergence. A relationship between the two  is provided by a theorem in the theory of $k$-Gevrey formal series, which  ensures that the formal $k$-Borel transform of $\tilde{\varphi} \in t\mathbb{C}\llb t \rrb$ is analytically convergent if and only if $\tilde{\varphi}$ is $k$-Gevrey. We can synthetically formulate  the theorem as \cite{aD15}:
\begin{dmath}\label{eq:k-gevrey-theorem}
	\op{B}_k \left( t\mathbb{C}\llb t \rrb_{\frac{1}{k}} \right) = \bbb{C}\{\zeta\}
\end{dmath}.

Note that, if $\varphi(t) \in \mathbb{C}\llb t \rrb_{\frac{1}{k}}$, then trivially $\tilde{\varphi}(t) = t\,\varphi(t) \in t\mathbb{C}\llb t \rrb_{\frac{1}{k}}$. 
So, theorem \labelcref{eq:k-gevrey-theorem} ensures that, for any $\varphi(t) \in \mathbb{C}\llb t \rrb_{\frac{1}{k}}$, the formal $k$-Borel transform of $\tilde{\varphi}(t) = t\,\varphi(t)$ belongs to $\bbb{C}\{\zeta\}$, i.e. $\tilde{\varphi}_{\op{B}_k}(\zeta)$ is analytically convergent. 

Suppose that $\tilde{\varphi}_{\op{B}_k}(\zeta)$ can be  extended by analytic continuation to a function (which, for simplicity, we identify with the same symbol), holomorphic in $V = S(d, \delta)$, an infinite open sector with vertex in the origin, bisected by $d$ and of aperture $\delta$. 

\noindent Further suppose that the extended function $\tilde{\varphi}_{\op{B}_k}(\zeta)$ grows to infinity at most exponentially along $d$.

\noindent Then, according to the theory of Borel-Laplace resummation, the $k$-Gevrey formal series $\varphi$ is $k$-summable in the direction $d$. Its $k$-sum is a function $\op{S}[\varphi](t)$, defined in terms of the $k$-Laplace transform\footnote{
		Note that, if $\tilde{\varphi} \in t\bbb{C}\{t\}$, then $\mathcal{L}_k[\op{B}_k [\tilde{\varphi}]] = \tilde{\varphi}$, that is the $k$-Laplace transform is the inverse of the $k$-Borel transform in the subalgebra of convergent formal series without constant term. \label{fn:borel-laplace}
		} of $\tilde{\varphi}_{\op{B}_k}(\zeta)$ in the direction $d$ by the formula \cite{wB83, eM03}:
\begin{dmath}\label{eq:k-Borel_sum_definition}
	\op{S}[\varphi](t) \defeq {
		\frac{1}{t}\mathcal{L}_k[\tilde{\varphi}_{\op{B}_k}](t)  
	}
\end{dmath}.
Furthermore, it is always possible to extend $\op{S}[\varphi]$ to an open sector $V = S(d, \Theta)$ of aperture $\Theta > \pi/k$ such that $\varphi \sim_{\overset{1/k}{V}} \op{S}[\varphi]$. This means that $\op{S}[\varphi]$ is the unique function $k$-asymptotic to $\varphi$ in the direction $d$, where $k$ is the exact Gevrey level of $\varphi$.

By definition, a $k$-Gevrey formal series is $k$-summable (in the sense of Borel-Laplace resummation) if it is $k$-summable in almost all directions, i.e. in all directions excepted at most a subset of measure zero. The subalgebra of $k$-summable formal series is denoted $\bbb{C}\{t\}_{\frac{1}{k}}$ and, for different $k, k'$ the following Tauberian property holds:
\begin{dmath*}
		\bbb{C}\{t\}_{\frac{1}{k}} \cap \bbb{C}\{t\}_{\frac{1}{k'}} = \bbb{C}\{t\}
\end{dmath*}.

\section{An umbral theory in $\mathbb{C}\llb t\rrb$}\label{sec:umbral_theory}

Before discussing how umbral concepts can be formalised within the context of formal power series, we will briefly review the lines along which the indicial method has developed so far (see \cite{sL22} for an exhaustive treatment).

The very idea of "umbral operator" stems from the observation that $\op{e}^{\mu \partial_z}$ formally produces a shift by $\mu$ in the argument of sufficiently smooth functions. In particular, the identity
\begin{dmath*}
	\op{e}^{\mu \partial_z} \varphi(z)\Big|_{z=0} = \varphi(\mu)
\end{dmath*}
is reinterpreted as the action of an informally defined "operator $\uo$ raised to the power $\mu$" on an equally informally defined "ground state" $\bvarphi$:
\begin{dmath}\label{eq:umbral_operator}
	\uo^{\mu} \bvarphi = \varphi(\mu)
\end{dmath}.

\noindent \Cref{eq:umbral_operator} is then exploited in expressions like:
\begin{dmath*}
	F(w) = {f(w\uo^\mu)\,\bvarphi }
\end{dmath*},
where both $F$ and $f$ are assumed to admit Maclaurin expansions.

In other terms, a heuristic identity is established between some highly transcendental function $F$ and a simpler function $f$\footnote{
	The function $f$ is the "umbra" left by $F$ under the action of $\uo$, thus justifying the name of the method. 
	}, which becomes manifest when the Maclaurin expansion of the latter is "evaluated" in $w \uo^\mu$ and the resulting operator applied to some ground state $\bvarphi$ identified \emph{ad hoc}.

\subsection{The umbral operator}

Any acceptable definition of umbral operator in the context of formal series must be able to produce a rigorous analog of \cref{eq:umbral_operator}.

One might be led to consider as a promising starting point the following composition\footnote{Note that the composition is defined as a series of formal series.
} of formal series:
\begin{dmath}\label{eq:formal_composition}
	\varphi\circ(\op{id} + \chi) \defeq \sum_{n=0}^\infty \frac{1}{n!}\,\chi^n \partial^n\varphi \condition*{\varphi$, $\chi \in \mathbb{C}\llb t \rrb}
\end{dmath},
which, when $\chi = \mu \in \bbb{C}$, is expected to produce the shift $\varphi(t) \mapsto \varphi(t+\mu)$.

Unfortunally, \cref{eq:formal_composition} is well defined -- i.e. converges formally in the Krull topology --  \emph{only if} $\op{val}(\chi) \geq 1$, due to the valuation-lowering nature of $\partial$. In particular, it is not formally convergent when $\chi = \mu \in \bbb{C}$, so a formal shift operator acting as
$(\op{T}_\alpha\varphi)(t) \defeq \varphi(t + \mu)$ cannot be generally defined in $\mathbb{C}\llb t \rrb$.

On the other hand, the "evaluation at zero" operator, $\uo^0 \,:\mathbb{C}\llb t \rrb \to \bbb{C} \subset\mathbb{C}\llb t \rrb$, can be rigorously defined as the  functional that returns the constant term of its argument:
\begin{dmath*}
	\varphi(t) \mapsto {\uo^0[\varphi] \equiv \varphi(0) \defeq c^{(\varphi)}_0}
\end{dmath*}.
By exploiting $\uo^0$, we can give a rigorous definition of the "evaluation at $\mu$" operator $\uo^{\mu}: \mathbb{C}\llb t \rrb \to \bbb{C}\subset\mathbb{C}\llb t \rrb$ as the functional:
\begin{dmath}\label{eq:evaluation_at_alpha}
	\varphi(t) \mapsto {\uo^{\mu}[\varphi] \defeq \sum_{n=0}^\infty \frac{\mu^n}{n!}\uo^0[\partial^n\varphi] \equiv \sum_{n=0}^\infty \frac{\mu^n}{n!}c^{(\varphi)}_n =\varphi(\mu)}
\end{dmath}.

The definition is well posed \emph{provided} $\varphi(\mu) \in \bbb{C}$. In other words, the formal series $\varphi$, interpreted as a power series in $\bbb{C}$, must converge in $\mu$ in the ordinary analytic sense,  $\varphi \in \bbb{C}\{t\}$. Assuming the convergence of $\varphi$  within an open disk of radius $R_\varphi$ (possibly infinite, if $\varphi$ is entire), then the definition \cref{eq:evaluation_at_alpha} is well posed for any $\mu$ such that\footnote{
		The definition can be possibly extended to $|\mu|\, > R_\varphi$ if $\varphi$ -- as it is often the case -- is analytically continuable to larger regions of the complex plane. Of course, any value of $\mu$ is in principle allowed when $\varphi$ is entire.
}  $|\mu|\, < R_\varphi$.

We call $\uo = \uo^1$ \emph{umbral operator} and $\varphi$ \emph{umbral ground state}, and require that the parameter $\mu$ behaves as an exponential index (whence the qualification "indicial" for the method). Namely:
\begin{dgroup*}
	\begin{dmath*}
		(\uo^{\mu} \uo^{\,\nu})[\varphi] = {\uo^{\mu + \nu} [\varphi] = \varphi(\mu + \nu)}
	\end{dmath*},
	\begin{dmath*}
		((\uo^{\mu})^{\,\nu})[\varphi] = {\uo^{\mu\nu} [\varphi] = \varphi(\mu\nu)}
	\end{dmath*},
\end{dgroup*}
provided the corresponding numeric series are convergent.

It is important to stress that this transposition of the indicial umbral method in the wider context of formal power series only requires the introduction of a \emph{single} umbral operator $\uo$. By acting on different umbral ground states $\varphi \in \bbb{C}\{t\}$ with this single operator, it is possible to reproduce all the case studies considered so far in the literature on the topic, where a plethora of \emph{ad hoc} different umbral operators are defined and used.

\subsection{Umbral identities}

In order to lay solid foundations for the indicial umbral method and fully unleash its potential, the next step is  to make rigorous sense of umbral identities of the type:
\begin{dmath}\label{eq:umbral_identity}
	F(\zeta) = {f(\zeta\uo^\mu)[\varphi] }
\end{dmath},
which are at the very core of the method and often appear in its applications.

Before giving a sound definition of \cref{eq:umbral_identity} in the context of of formal power series, we make some preliminary assumptions:
\begin{enumerate}
	\item $\uo: \mathbb{C}\llb t \rrb \to \bbb{C}\subset\mathbb{C}\llb t \rrb$ is the umbral operator defined by \cref{eq:evaluation_at_alpha};
	\item $\mu \geq 0$;
	\item the ground state $\varphi \in \bbb{C}\{t\}$ is a formal series analytically convergent in  $D(R_\varphi)$, an open disk centred in the origin of (possibly infinite) radius $R_\varphi$;
	\item $f \in \bbb{C}\{t\}$ is also an analytically convergent formal series\footnote{This assumption is not strictly necessary. It is introduced here in order to make contact with the concrete applications of the umbral framework, where $f$ is normally a function.}, expressible in terms of the Maclaurin expansion:
		\begin{dmath*}
		f(t) = {
			\sum_{r = 0}^\infty c^{(f)}_r t^r =
			\sum_{r=0}^\infty \frac{f^{(r)}(0)}{r!}\,t^{r} \condition*{f^{(r)}(0) \defeq (\partial_t^r f)(0)}
		}
	\end{dmath*},
	which can be
	\begin{enumerate}
		\item convergent in a disk $|t|< R$; 
		\item convergent for any $t \in \bbb{C}$. In this case, we further assume that $f$ is of exponential type $A$, i.e. real constants $C, A$ exist such that $|f(t)| \leq C \op{e}^{A|t|}$ for $|t| \to \infty$.
	\end{enumerate}
\end{enumerate}
With these assumptions, we pose by definition:
\begin{dmath}\label{eq:rigorous_umbral_image}
	f(\zeta\uo^\mu)[\varphi] \defeq {
		\sum_{r = 0}^\infty c^{(f)}_r \uo^{\mu r}[\varphi]\,\zeta^r =
		\sum_{r=0}^\infty \frac{f^{(r)}(0)}{r!}\uo^{\mu r}[\varphi]\,\zeta^r =
		F(\zeta),
	}
\end{dmath}
for any $\mu \geq 0$\footnote{When $\mu = 0$, we have trivially $f(\zeta\uo^0)[\varphi] = f(\zeta)\,\varphi(0)$.} such that the numeric series $\uo^{\mu r}[\varphi]$ converges to $\varphi(\mu r) \in \bbb{C}$ for each $r \in \bbb{N}^0$. For such values of $\mu$, $F = f(\zeta\uo^\mu)[\varphi] \in \mathbb{C}\llb \zeta \rrb$ is a well-formed formal series in the indeterminate $\zeta$.

\noindent In other words, $F(\zeta)$ is the formal series resulting from the application to $\varphi$ of the operator $f(\zeta\uo^\mu): \bbb{C}\{t\} \to \bbb{C}\llb \zeta \rrb$, defined in terms of the functional $\uo^{\,\mu}: \bbb{C}\{t\} \to \bbb{C}$.

Note that, as a formal series, \cref{eq:rigorous_umbral_image} can be manipulated term by term with no concerns related to analytic convergence, provided the result of the manipulation is still a well-formed formal series.  

\noindent In particular, \cref{eq:rigorous_umbral_image} can be derived and integrated term by term with respect to $\zeta$. Since  by hypothesis $f$ is convergent to a function, we can derive or integrate directly this function, by treating $\uo^\mu$ as a parameter, provided we interpret the result formally, i.e. as a new formal series. 

\noindent This is a crucial property in applications of the umbral method, justified by the fact that the resulting formal series coincides with the one that is obtained by operating directly on \cref{eq:rigorous_umbral_image} term by term. This property is widely used to simplify calculations involving $F$, that can be carried out by replacing $F$ with the typically much simpler function $f$.

\begin{example}
As a trivial example, consider the formal series
\begin{dmath*}
	f(t) = {
		\sum_{r=0}^\infty \frac{c^r}{r!}\,t^r
	}
\end{dmath*},
converging to the entire function $\op{e}^{ct}$.
In this case \cref{eq:rigorous_umbral_image} reads:
\begin{dmath}\label{eq:example_exp}
	f(\zeta\uo^\mu)[\varphi]= {
		\op{e}^{c\,\zeta\uo^\mu}[\varphi]  =
		\sum_{r=0}^\infty \frac{c^r \uo^{\mu r}[\varphi]}{r!}\,\zeta^r
	}
\end{dmath}.
The derivative term by term of \cref{eq:example_exp} is:
\begin{dmath*}
	\partial \sum_{r=0}^\infty \frac{c^r \uo^{\mu r}[\varphi]}{r!}\,\zeta^r = 
	{
		\sum_{r=1}^\infty \frac{c^r \uo^{\mu r}[\varphi]}{r!}\,r\,\zeta^{r-1} =
		c\uo^\mu \op{e}^{c\,\zeta\uo^\mu}[\varphi]
	}
\end{dmath*},
which coincides with result obtained by directly deriving $ \op{e}^{c\,\zeta\uo^\mu}[\varphi]$, with $\uo^\mu$ treated as a parameter.

\end{example}
\vskip 0.2cm

\section{Umbral ground states}\label{sec:umbral_ground_states}

The definition of umbral identities as formal series, even if theoretically sound, is of little practical use if the series do not converge analytically. So we are naturally led to investigating the convergence conditions of \cref{eq:rigorous_umbral_image}. These obviously depend on the behaviour of the coefficients $f^{(r)}(0)\uo^{\mu r}[\varphi]/r!$, that is, ultimately, on the value of $\mu$ and the nature of the ground state $\varphi$.

For obvious practical reasons and in view of establishing a connection between the umbral formalism outlined in the previous section and the theory of Borel-Laplace resummation of formal series (see \cref{sec:borel_laplace_resummation}), we will limit for the moment the analysis to the following classes of ground state functions:
\begin{enumerate}
	\item the class of entire functions:
		\begin{dmath*}
			\phi_{\alpha,\, \beta}(t) \defeq {
			\frac{1}{\Gamma(\alpha t + \beta)} 
			\condition*{\alpha >0,\; \beta \in \bbb{C}}
			}
		\end{dmath*};
	\item the class of meromorphic functions:
		\begin{dmath*}
		\psi_{\alpha,\,\beta\,,\gamma}(t) = {
			\frac{(\gamma)_t}{\Gamma(\alpha t + \beta)} \equiv
			\frac{\Gamma(\gamma + t)}{\Gamma(\gamma)\,\Gamma(\alpha t + \beta)}
			\condition*{\alpha >0,\; \beta, \gamma \in \bbb{C}}
		}
		\end{dmath*},
		where $(\gamma)_t \defeq \Gamma(\gamma + t)/\Gamma(\gamma)$ is a generalised version of the rising Pochhammer symbol. $\psi_{\alpha,\,\beta\,,\gamma}(t)$  has simple poles when $\gamma + t \in \{0\} \cup \bbb{Z}^-$.
\end{enumerate}

It is perhaps worth remarking that restricting the analysis to these functions is not a severe limitation, since the vast majority of ground states considered so far in the literature on the topic are of this type.

In the following subsections we analyse the convergence conditions of \cref{eq:rigorous_umbral_image} determined by these classes of ground state functions. Note that both classes are characterised by the presence of the Gamma function\footnote{See e.g. \cite{gD10} for a review of the pivotal role played by the Euler Gamma function in the development of many aspects of modern physics.} $\Gamma(\alpha t + \beta)$ at the denominator, a crucial factor for establishing a connection with the theory of Borel-Laplace resummation.

\subsection{The ground state $\phi_{\alpha,\, \beta}(t)$}

With this specific ground state, \cref{eq:rigorous_umbral_image} becomes:
\begin{dmath}\label{eq:rigorous_umbral_image_phi_alpha_beta}
	F(\zeta) = {
		f(\zeta\uo^\mu)[\phi_{\alpha,\, \beta}] =
		\sum_{r = 0}^\infty \frac{c^{(f)}_r}{\Gamma(\mu\alpha r + \beta)} \,\zeta^r
		\condition*{\alpha > 0,\, \beta \in \bbb{C}}
	}
\end{dmath},
where $c^{(f)}_r = f^{(r)}(0)/r!$. Note that, although the reciprocal of the Gamma function is well defined for any value of the argument, we exclude combinations of the parameteres such that $\mu\alpha r +\beta \in \{0\} \cup \bbb{Z}^-$, which correspond to values of $r$ where the coefficients of $F$ are formally undefined . 

\subsubsection{Asymptotic analysis of the case  $\alpha = \beta=1$}
We first consider the special case $\alpha = \beta=1$, that is following the archetypical ground state function, often appearing in applications:
\begin{dmath}\label{eq:umbral_state_phi}
	\phi(t) \defeq {\phi_{1,1}(t) = \frac{1}{\Gamma(1+t)}}
\end{dmath}.
The asymptotic analysis of the convergence conditions of the resulting formal series:
\begin{dmath}\label{eq:rigorous_umbral_image_phi}
	F(\zeta) = {
		f(\zeta\uo^\mu)[\phi] =
		\sum_{r = 0}^\infty \frac{c^{(f)}_r}{\Gamma(1+\mu r)} \,\zeta^r
	}
\end{dmath},
is  carried out in detail, by treating separately the cases when $f$ is simply convergent or entire of exponential type.
\begin{enumerate}
	\item If $f$ is holomorphic in a disk $|\,t\,| < R$, by Cauchy's estimates there exists \( C > 0 \) such that $| c^{(f)}_r |\,\leq C/R^r$. Then
		\begin{dmath*}
			|c^{(F)}_r| = {
				\left| \frac{c^{(f)}_r}{\Gamma(1 + \mu r)} \right| \leq \frac{C}{R^r \,\Gamma(1 + \mu r)}.
			}
		\end{dmath*}
		Since $\Gamma(1 + \mu r) \to \infty$ super-exponentially, $c^{(F)}_r \to 0$ for $r \to \infty$. So the series $F(\zeta)$ converges for all $\zeta \in \mathbb{C}$ -- hence it is of 0-Gevrey.
	\item if $f$ is entire of exponential type $A>0$, then for any $r$
		\begin{dmath*}
			|c^{(f)}_r| \leq C A^r \frac{\op{e}^r}{r^r} \Rightarrow  |c^{(F)}_r| \,\leq C\,\frac{A^r\op{e}^r/r^r}{\Gamma(1 + \mu r)} 
		\end{dmath*}. 
		Using Stirling's formula:
		\begin{dmath*}
			\Gamma(1 + \mu r) \sim \sqrt{2\pi \mu r} \cdot (\mu r)^{\mu r} e^{-\mu r}\condition*{r \to \infty}
		\end{dmath*},
		we obtain the asymptotic estimate
		\begin{dmath*}
			|c^{(F)}_r|\, \lesssim \left( \frac{K}{r^{1+\mu}} \right)^r \quad \text{for some constant } K > 0
		\end{dmath*}.
		Also in this case the series converges absolutely for all $\zeta \in \mathbb{C}$ and any $\mu>0$.
\end{enumerate}
In conclusion, when $\varphi = \phi$, \cref{eq:rigorous_umbral_image_phi} converges to an entire function $F(\zeta)$ for any $\mu > 0$. 

\subsubsection{Gevrey analysis of the case  $\alpha = \beta=1$}
Remarkably, the same conclusion can be reached by directly resorting to the Gevrey theory sketchily outlined in the previous section.

\noindent In fact, by direct comparison with \cref{eq:rigorous_borel_transform_k}, we notice that $F$ in \cref{eq:rigorous_umbral_image_phi} is the $\frac{1}{\mu}$-Borel transform of the 0-Gevrey formal series $\tilde{f}(t) = t \,f(t)$. Because of the nesting property \labelcref{eq:gevrey_nested}, $\tilde{f}$ is also $\frac{1}{\mu}$-Gevrey -- i.e. $\tilde{f} \in t\bbb{C}\llb t \rrb_\mu$ -- for any $\mu > 0$. Then theorem \labelcref{eq:k-gevrey-theorem} guarantees that its $\frac{1}{\mu}$-Borel transform $F$ is convergent for any $\mu > 0$. 

\noindent $F$ is in fact not only simply convergent, but even entire (as directly deducible from the asymptotic analysis), because the growth of the coefficients of the $\frac{1}{\mu}$-Borel transform of a formal series $\tilde{f} \in t\bbb{C}\llb t \rrb_0$ is strictly subfactorial.

\subsubsection{General case  $\alpha>0$, \;$\beta \in \bbb{C}$}
The previous observations enable to analyse with minimal effort the convergence conditions of \cref{eq:rigorous_umbral_image_phi_alpha_beta} for $\mu > 0$ and general values of the parameters $\alpha > 0,\, \beta \in \bbb{C}$, with $\mu\alpha r +\beta \notin \{0\} \cup \bbb{Z}^-$.

\noindent Indeed, by direct comparison with \cref{eq:rigorous_borel_transform_k_beta},  it is straightforward to recognise that \cref{eq:rigorous_umbral_image_phi_alpha_beta} is the generalised $(\frac{1}{\alpha\mu}, \beta)$-Borel transform of the already introduced 0-Gevrey formal series $\tilde{f} \in t\bbb{C}\llb t \rrb_0$. Since the presence of the shift $\beta$ does not change the Gevrey level of the transform, by the same argument used before $F$ converges to en entire function for any $\mu> 0, \,\alpha > 0$ and $\beta \in \mathbb{C} \setminus \{-\alpha\mu r : r \in \mathbb{N}^0\}$.

\begin{example}

Let  $f(t)$ be the entire function $\op{e}^{-t}$. Using this $f$ and the ground state $\phi$, and setting $\mu=1$, \cref{eq:rigorous_umbral_image_phi} reads:
\begin{dmath*}
	F(\zeta) = {
		\sum_{r=0}^\infty \frac{(-1)^r}{(r!)^2}\,\zeta^r
	}
\end{dmath*},
which coincides with the 1-Borel transform of $\tilde{f}(t) = t\op{e}^{-t}$. This power series converges for any $\zeta \in \bbb{C}$ to $C_0(\zeta)$, the Tricomi function of order 0. 

Note that the ordinary 1-Laplace transform of this function yields as expected (see footnote \ref{fn:borel-laplace}):
\begin{dmath*}
	\mathcal{L}_1[C_0](t) = t\op{e}^{-t} \condition*{\op{Re}t >0}
\end{dmath*},
which can be continued to the entire complex plane.
So, the following umbral identity holds for any $\zeta \in \bbb{C}$:
\begin{dmath*}
	C_0(\zeta) = \op{e}^{-\zeta\uo}[\phi]
\end{dmath*}.

\end{example}
\vskip 0.2cm

\begin{example}
Let $f(t) = 1/(1 - t^2)$, whose Maclaurin power expansion converges in a disk $|\,t\,| \,< 1$.
Using this $f$ and the ground state $\phi$, \cref{eq:rigorous_umbral_image_phi} with $\mu=1$ reads:
\begin{dmath*}
	F(\zeta) = {
		\sum_{r=0}^\infty \frac{z^{2r}}{(2r)!}
	}
\end{dmath*},
which coincides with the 1-Borel transform of $\tilde{f}(t) = t/(1-t^2)$ and the Maclaurin of expansion of $\cosh(\zeta)$, an entire function.

\noindent Also in this case, as expected:
\begin{dmath*}
	\mathcal{L}_1[\cosh](t) = \frac{t}{1 - t^2} \condition*{0< \op{Re}t < 1, \,\op{Im} t = 0}
\end{dmath*},
which can be continued to $|\,t\,| \,< 1$.
So, the following umbral identity holds for any $\zeta \in \bbb{C}$:
\begin{dmath*}
	\cosh(\zeta) = \frac{1}{1 - (\zeta\uo)^2}[\phi] \end{dmath*}.
\end{example}
\vskip 0.2cm

\begin{example}

Let  $f(t)$ be the entire function $\op{e}^{t}$. Using this $f$ and the general ground state $\phi_{\alpha,\, \beta}$, \cref{eq:rigorous_umbral_image_phi} with $\mu=1$ reads:
\begin{dmath*}
	F(\zeta) = \sum_{r = 0}^\infty \frac{z^r}{r! \Gamma(\alpha r + \beta)} \condition*{\alpha >0,\; \beta \in \bbb{C}}
\end{dmath*},
which coincides with the $(\frac{1}{\alpha},\,\beta)$-Borel transform of $\tilde{f}(t) = t\op{e}^{t}$. This power series converges for any $\zeta \in \bbb{C}$, and is in fact the definition of the Bessel-Wright function $W_{\alpha,\,\beta}(\zeta)$ \cite{emW34}, an entire function of $\zeta$.

Once again, the $(\frac{1}{\alpha},\,\beta)$-Laplace transform of this function yields as expected:
\begin{dmath*}
	\mathcal{L}_{\frac{1}{\alpha},\,\beta}[W_{\alpha,\,\beta}](t) = t\op{e}^{t} 
\end{dmath*}.
So, the following umbral identity holds for any $\zeta \in \bbb{C}$:
\begin{dmath*}
	W_{\alpha,\,\beta}(\zeta) = \op{e}^{\zeta\uo}[\phi_{\alpha,\, \beta}]
\end{dmath*}.

\end{example}
\vskip 0.2cm

\subsection{The ground state $\psi_{\alpha,\,\beta,\,\gamma}(t)$}

We continue our analysis with the class of ground states $\psi_{\alpha,\,\beta,\,\gamma}$, by considering first the special case $\gamma = 1$.

\subsubsection{Case $\gamma=1$}

Since $(1)_t = \Gamma(1 + t)$, the choice $\gamma = 1$ corresponds to the the ground state:
\begin{dmath*}
	\psi_{\alpha,\,\beta}(t) \defeq {
	\psi_{\alpha,\,\beta,\,1}(t) =
		\frac{\Gamma(1 + t)}{\Gamma(\alpha t + \beta)}
	}
\end{dmath*},
a meromorphic function with simple poles when $t = - n$, \, $n \in \bbb{N}$. 

Using this ground state, \cref{eq:rigorous_umbral_image} becomes:
\begin{dmath}\label{eq:rigorous_umbral_image_psi_alpha_beta}
	F(\zeta) = {
		f(\zeta\uo^\mu)[\psi_{\alpha,\,\beta}] =
		\sum_{r = 0}^\infty c^{(f)}_r\frac{\Gamma(\mu r + 1)}{\Gamma(\alpha\mu r + \beta)} \,\zeta^r
	}
\end{dmath},
where $c^{(f)}_r = f^{(r)}(0)/r!$.

\noindent Provided the additional condition $\alpha\mu r + \beta \notin \{0\}\cup\bbb{Z}^-$ is satisfied for all $r \in \bbb{N}^0$, \cref{eq:rigorous_umbral_image_psi_alpha_beta} is recognised as the  $(\frac{1}{\alpha\mu},\,\beta)$-Borel transform of the formal series $\tilde{g}(t) = t\,g(t) \in t\bbb{C}\llb t\rrb$, where:
\begin{dmath*}
	g(t) = {
	\sum_{r = 0}^\infty c^{(g)}_r t^r \defeq
	\sum_{r = 0}^\infty c^{(f)}_r\,\Gamma(\mu r + 1)\, t^r
	}
\end{dmath*}.

The asymptotic analysis of the coefficients $c^{(g)}_r$, carried out by approximating the Gamma function with the Stirling formula, shows that:
\begin{enumerate}
	\item $g \in \bbb{C}\llb t\rrb_0$, i.e. $g$ is analytically convergent, if $f$ is analytically convergent in a disk of finite radius;
	\item $g \in \bbb{C}\llb t\rrb_\mu$, i.e. $g$ is divergent of Gevrey level $\frac{1}{\mu}$, if $f$ is entire.
\end{enumerate}
\noindent In virtue of property \labelcref{eq:gevrey_nested}, in both cases $\tilde{g}$ is also \emph{a fortiori} $\frac{1}{\alpha\mu}$-Gevrey if $\alpha \geq 1$. Then theorem \labelcref{eq:k-gevrey-theorem} ensures that its $(\frac{1}{\alpha\mu},\,\beta)$-Borel transform -- which coincides with $F$ -- is always convergent for $\alpha \geq 1$.

\noindent This result is confirmed by the direct asymptotic analysis of the coefficients $c^{(F)}_r$, which demonstrates that \cref{eq:rigorous_umbral_image_psi_alpha_beta} is convergent to an entire function when $\alpha \geq 1$. The term $\Gamma(\alpha\mu r + \beta)$ prevails over the coefficients of $g$, which diverge in modulus at most as $(r!)^{\frac{1}{\mu}}$, thus ensuring an exponential convergence.

On the other hand, when $\alpha < 1$:
\begin{enumerate}
	\item if $g \in \bbb{C}\llb t\rrb_0$ then $F$ converges to an entire function;
	\item if $g \in \bbb{C}\llb t\rrb_\mu$, then $F$ is a divergent formal series of Gevrey level $\frac{1}{\mu(1-\alpha)}$, that is $F \in \bbb{C}\llb \zeta\rrb_{\mu(1-\alpha)}$.  
\end{enumerate}

\noindent Indeed, if $g \in \bbb{C}\llb t\rrb_\mu$ -- hence $\tilde{g} \in t\bbb{C}\llb t\rrb_\mu$ -- the Gamma function at the denominator does not fully suppress the factorial growth of the coefficients of $g$\footnote{
	Intuitively, since $|c^{(\tilde{g})}_r| \lesssim  (r!)^{\,\mu}$ and $\Gamma(1+\alpha\mu r)\lesssim  (r!)^{\alpha\mu}$, then $|c^{(F)}_r| \lesssim  (r!)^{\,\mu(1-\alpha)}$.
	}.  In this case, assuming that all the accessory conditions are fulfilled, we can resort to the theory of Borel-Laplace resummation and find, if it exists, the unique function $k$-asymptotic to $F$ in an open sector $S(d, \Theta)$ of aperture $\Theta > \pi/k$ and bisected by $d$ -- typically, the positive real line -- with $k = 1/\mu(1-\alpha)$.

In conclusion:
\begin{enumerate}
	\item when $\alpha \geq 1$, $F$ is an entire function of $\zeta$, hence 0-Gevrey;
	\item when $\alpha < 1$, $F$ either converges to an entire function if $g \in \bbb{C}\llb t\rrb_0$ or is divergent of Gevrey level $\frac{1}{\mu(1-\alpha)}$ if $g \in \bbb{C}\llb t\rrb_\mu$. In the latter case, it possibly represents the $\frac{1}{\mu(1-\alpha)}$-asymptotic expansion of a function that can be found by means of Borel-Laplace resummation.
\end{enumerate}
Some examples will help get a better feeling of the whole process.

\begin{example}

Let  $f(t)$ be the entire function $\op{e}^{t}$. Assuming $\mu=1$ and using this $f$ and the ground state $\psi_{\alpha,\, \beta}$, \cref{eq:rigorous_umbral_image_psi_alpha_beta} becomes:
\begin{dmath}\label{eq:ML_series}
	F(\zeta) = {
		\op{e}^{\zeta\uo} \,[\psi_{\alpha,\, \beta}] =
		\sum_{r = 0}^\infty \frac{z^r}{\Gamma(\alpha r + \beta)} \condition*{\alpha >0,\; \beta \in \bbb{C}}
	}
\end{dmath},
which coincides with the $(\frac{1}{\alpha},\,\beta)$-Borel transform of $\tilde{g}(t) = t\,g(t) = t\sum_{r = 0}^\infty t^r$, converging to $t/(1-t)$ for $|\,t\,|\,<1$.

Since $\tilde{g} \in t\bbb{C}\llb t\rrb_0$, property \labelcref{eq:gevrey_nested} and theorem \labelcref{eq:k-gevrey-theorem} ensure that $F$ converges for any $\zeta \in \bbb{C}$ and any $\alpha > 0$. \Cref{eq:ML_series} is in fact the Maclaurin expansion of the two-parameters Mittag-Leffler function $E_{\alpha,\,\beta}(\zeta)$, an entire function of $\zeta$ (see e.g. \cite{tN23} for a conventional umbral approach).

The $(\frac{1}{\alpha},\,\beta)$-Laplace transform of this function yields as expected:
\begin{dmath*}
	\mathcal{L}_{\frac{1}{\alpha},\,\beta}[E_{\alpha,\,\beta}](t) = \frac{t}{1-t} \condition*{|\,t\,|\,<1}
\end{dmath*}.

In conclusion, the following umbral identity holds for any $\zeta \in \bbb{C}$:
\begin{dmath*}
	E_{\alpha,\,\beta}(\zeta) = \op{e}^{\zeta\uo}[\psi_{\alpha,\, \beta}]
\end{dmath*}.

\end{example}
\vskip 0.2cm

\begin{example}

Let  $f(t)$ be the entire function $\op{e}^{it}$. Assuming $\mu=1$ and using this $f$ and the ground state 
\begin{dmath*}
	\lambda \defeq {
		\psi_{\frac{1}{2},\, 1} =
		\frac{\Gamma(1+t)}{\Gamma(1 +t/2)}
	}
\end{dmath*},
\cref{eq:rigorous_umbral_image_psi_alpha_beta} becomes:
\begin{dmath}\label{eq:faddeeva_series}
	F(\zeta) = {
		\op{e}^{i\zeta\uo} \,[\lambda] =
		\sum_{r = 0}^\infty \frac{i^{\,r}} {\Gamma(1 + r/2)}\, \zeta^r 
	}
\end{dmath}.
This represents an explicit example of \cref{eq:rigorous_umbral_image_psi_alpha_beta} with $\alpha = 1/2 < 1$. The formal series $F$ coincides with the $2$-Borel transform of $\tilde{g}(t) = t\,g(t) = t\sum_{r = 0}^\infty (it)^r$, converging to $t/(1-it)$ for $|\,t\,|\,<1$. 

Since $\tilde{g} \in t\bbb{C}\llb t\rrb_0$, property \labelcref{eq:gevrey_nested} and theorem \labelcref{eq:k-gevrey-theorem} ensure also in this case that $F$ converges to an entire function of $\zeta$.
\Cref{eq:faddeeva_series} is in fact the Maclaurin expansion of the Faddeeva function  \cite{vnF54} 
\begin{dmath*}
	\op{w}(\zeta) = {
		E_{\frac{1}{2},\, 1}(i\zeta) =
		\op{e}^{-\zeta^2}(1+i\,\op{erfi}(\zeta))
	}
\end{dmath*},
 also known as Gaussian exponential function (see \cref{sec:gaussian_trig} below).

The 2-Laplace transform of this function yields as expected:
\begin{dmath*}
	\mathcal{L}_2[\op{w}](t) = \frac{t}{1-it} \condition*{|\,t\,|\,<1}
\end{dmath*}.
In conclusion, the following umbral identity holds for any $\zeta \in \bbb{C}$:
\begin{dmath}\label{eq:faddeeva_function}
	\op{w}(\zeta) = \op{e}^{i\zeta\uo}[\lambda]
\end{dmath}.

\end{example}
\vskip 0.2cm

\begin{example}\label{ex:7}

Let us consider another instructive example of \cref{eq:rigorous_umbral_image_psi_alpha_beta} with $\alpha = 1/2 < 1$. 

Take $f(t) = 1/(1 +t)$, whose Maclaurin expansion converges in $|\,t\,| < 1$. Assuming $\mu=2$ and using this $f$ and the already introduced ground state $\lambda \defeq \psi_{\frac{1}{2},\, 1}$, \cref{eq:rigorous_umbral_image_psi_alpha_beta} becomes:
\begin{dmath}\label{eq:divergent_series}
	F(\zeta) = {
		\frac{1}{1 + \zeta\uo^2} \,[\lambda] =
		\sum_{r=0}^\infty(-1)^r\frac{\Gamma(2r+1)}{\Gamma(r+1)}\,\zeta^r =
		\sum_{r=0}^\infty(-1)^r\frac{(2r)!}{r!}\,\zeta^r
	}
\end{dmath}.
 \Cref{eq:divergent_series} is the $1$-Borel transform of $\tilde{g}(t) = t\,g(t) = t\sum_{r = 0}^\infty (2r)!(-t)^r$, a divergent series, of Gevrey level $ k = \frac{1}{\mu} = \frac{1}{2}$, that is $\tilde{g} \in t\llb \bbb{C}(t) \rrb_2$. 
 
 $F$ is then divergent of Gevrey level $\frac{1}{\mu(1-\alpha)}=1$, and one can try to make sense of it by means of the Borel-Laplace resummation procedure. We first take the 1-Borel transform\footnote{Recall that the Borel transform is only defined for formal series without constant term.} of $\tilde{F}(z) = z\,F(z)$, where we have changed the symbol for the indeterminate for later convenience. Namely:
 \begin{dmath*}
	\hat{F}_1(z) \defeq {
		\op{B}_1 [\tilde{F}](z) =
		\sum_{r = 0}^\infty \frac{(-1)^r(2r)!}{r!\Gamma(r + 1)}\,z^r =
		\sum_{r = 0}^\infty \frac{(-1)^r}{(r!)^2}\,z^r =
		(1+4z)^{-\frac{1}{2}}
	}
\end{dmath*}.
$\hat{F}_1$ converges to a function that, though not entire, can be analytically continued so as to satisfy all the  conditions needed for the existence of the 1-Laplace integral transform along $\bbb{R}^+$. By applying such transform (second step of the Borel-Laplace resummation procedure) we eventually obtain:
\begin{dmath*}
	\op{S}[F](\zeta) \defeq {
		\frac{1}{\zeta}\,\mathcal{L}_1[\hat{F}_1](\zeta) = 
		\frac{1}{\zeta} \int_0^\infty \mathrm{d}z\,(1+4z)^{-\frac{1}{2}} \op{e}^{-z/\zeta}
		} = 
		\frac{1}{2}\sqrt{\frac{\pi}{\zeta}}  \op{e}^{1/4\zeta} \op{erfc}(1/2\zeta^{1/2})\condition*{\zeta \in V=S(\Theta),\;\Theta=2\pi}
\end{dmath*}.

Note that $F \in \llb \bbb{C}(t) \rrb_1$ implies $F \in \llb \bbb{C}(t) \rrb_k$ for any $k \geq 1$, so the resummation can also in principle be performed by taking Borel and Laplace tranforms of order $\frac{1}{k}$, with $k \geq 1$. For example, for $k=2$ we obtain:
\begin{dmath*}
	\hat{F}_{\frac{1}{2}}(z) \defeq {
		\op{B}_{\frac{1}{2}} [\tilde{F}](z) =
		\sum_{r = 0}^\infty \frac{(-1)^r(2r)!}{r!\Gamma(2r + 1)}\,z^r =
		\sum_{r = 0}^\infty \frac{(-1)^r}{r!}\,z^r
	}
\end{dmath*}.
$\hat{F}_{\frac{1}{2}}$ converges to the entire function $\op{e}^{-z}$, which satisfies all the  conditions needed for the convergence of the $\frac{1}{2}$-Laplace transform along $\bbb{R}^+$. By applying such transform we eventually obtain:
\begin{dmath*}
	\op{S}[F](\zeta) \defeq {
		\frac{1}{\zeta}\,\mathcal{L}_{\frac{1}{2}}[\hat{F}](\zeta) = 
		\frac{1}{2\sqrt{\zeta}} \int_0^\infty \mathrm{d}z\op{e}^{-z} \,z^{-1/2} \op{e}^{-(z/\zeta)^{1/2}}} = 
		\frac{1}{2}{\sqrt{\frac{\pi}{\zeta}}} \op{e}^{1/4\zeta} \op{erfc}(1/2\zeta^{1/2})\condition*{\zeta \in V=S(\Theta),\;\Theta=2\pi}
\end{dmath*}.
The two procedures give the same result, as expected. Since $F$ is exactly 1-Gevrey, we say that $S[F](\zeta)$ is the unique function 1-Gevrey asymptotic to $F$ in $V$\footnote{Note that the aperture of $V$ satisfies the condition $\Theta = 2\pi > \pi/k = \pi$.}. Equivalently, the formal series $F(\zeta) = 1/(1 + \zeta \uo^2)\,[\lambda]$ represents the unique 1-asymptotic expansion of $S[F](\zeta)$ for $\zeta \in V$.

In conclusion, the following asymptotic umbral identity holds for $\zeta \in V$:
\begin{dmath*}
	\op{S}[F](\zeta) \sim_{\overset{1}{V}} \frac{1}{1 + \zeta\uo^2} \,[\lambda] 
\end{dmath*}.

\end{example}
\vskip 0.2cm

\begin{example}\label{ex:8}

Let us consider a further example, characterised by the same function $f$, ground state $\lambda$ and $\alpha = 1/2$ used in \cref{ex:7}, but a different value of $\mu$, namely $\mu=1$ instead of $\mu=2$.

In this case, \cref{eq:rigorous_umbral_image_psi_alpha_beta} becomes:
\begin{dmath}\label{eq:divergent_series_2}
	F(\zeta) = {
		\frac{1}{1 + \zeta\uo} \,[\lambda] =
		\sum_{r=0}^\infty(-1)^r\frac{\Gamma(r+1)}{\Gamma(1+r/2)}\,\zeta^r =
		\sum_{r=0}^\infty(-1)^r\frac{r!}{\Gamma(1+r/2)}\,\zeta^r
	}
\end{dmath}.
 \Cref{eq:divergent_series_2} is the $2$-Borel transform of $\tilde{g}(t) = t\,g(t) = t\sum_{r = 0}^\infty r!(-t)^r$, a divergent series of Gevrey level $\frac{1}{\mu} = 1$. 
 
The order $k=2$ of the Borel transform in \cref{eq:divergent_series_2} falls outside the domain of applicability of theorem \labelcref{eq:k-gevrey-theorem}, which would only ensure convergence if $0<k<1$. 
In fact, $F$ is divergent of Gevrey level $\frac{1}{\mu(1-\alpha)}=2$ , as can be directly verified by performing a thorough asymptotic analysis. 

The 2-Borel transform of $\tilde{F}(\zeta) = \zeta\,F(\zeta)$ is:
\begin{dmath*}
	\hat{F}(z) \defeq {
		\op{B}_{2} [\tilde{F}](z) =
		\sum_{r = 0}^\infty \frac{(-1)^r r!}{\Gamma(1+r/2)^2}\,z^r
	}
\end{dmath*},
converging to the function:
\begin{dmath*}
	\hat{F}(z) = \frac{2 \arccos(2 z)}{\pi  \sqrt{1-4 z^2}}
\end{dmath*}.
This function possesses branch points at $z = \pm 1/2$. For $t > 0$, the 2-Laplace integral:
\begin{dmath*}
    \mathcal{L}_{2}[\hat{F}](t) = \frac{2}{t} \int_0^{\infty} \zeta \hat{F}_1(\zeta) e^{-(\zeta/t)^2} d\zeta
\end{dmath*}
is ill-defined on the real axis. One must choose a different direction $d$ that avoids $z = 1/2$, leading to two distinct sectorial sums $S_{0^+}$ and $S_{0^-}$ which differ by an exponentially small factor:
\begin{dmath*}
    \Delta S \sim \op{e}^{-(1/2t)^2}
\end{dmath*}.
This is a manifestation of the Stokes phenomenon \cite{dS16}. We reserve to investigate this and other related aspects in a following paper.

\end{example}
\vskip 0.2cm

\subsubsection{General case}

For generic values of the parameter $\gamma$, \cref{eq:rigorous_umbral_image} becomes:
\begin{dmath}\label{eq:rigorous_umbral_image_phi_alpha_beta_gamma}
	F(\zeta) = {
		f(\zeta\uo^\mu)[\psi_{\alpha,\, \beta,\, \gamma}] = \frac{1}{\Gamma(\gamma)}
		\sum_{r = 0}^\infty \frac{c^{(f)}_r\, \Gamma(\mu r +\gamma)}{\Gamma(\mu\alpha r + \beta)} \,\zeta^r
		\condition*{\alpha > 0,\; \beta, \gamma \in \bbb{C}}
	}
\end{dmath}.
This is the $(\frac{1}{\alpha\mu},\, \beta)$-Borel transform of the formal series $\tilde{h}(t) = t\,h(t)$, where:
\begin{dmath}\label{eq:h_series}
	h(t) = {
		\sum_{r = 0}^\infty c^{(h)}_r\,t^r \defeq
		\sum_{r = 0}^\infty c^{(f)}_r \,\Gamma(\mu r +\gamma)\,t^r
	}
\end{dmath}.
The presence of the shift $\gamma$ in \cref{eq:h_series} has no influence on its Gevrey level, which is still either 0 or $\frac{1}{\mu}$, as in the particular case analysed in the previous subsection.

So, we can conclude with no futher ado that, also in the general case of \cref{eq:rigorous_umbral_image_phi_alpha_beta_gamma}:
\begin{enumerate}
	\item when $\alpha \geq 1$, $F$ converges to an entire function of $\zeta$, hence is 0-Gevrey;
	\item when $\alpha < 1$, $F$ either converges to an entire function if $g \in \bbb{C}\llb t\rrb_0$ or is divergent of Gevrey level $\frac{1}{\mu(1-\alpha)}$ if $g \in \bbb{C}\llb t\rrb_\mu$. In the latter case, it possibly represents the $\frac{1}{\mu(1-\alpha)}$-asymptotic expansion of a function that can be found by means of Borel-Laplace resummation.
\end{enumerate}
Of course, even if the shift $\gamma$ does not change the convergence conditions of $F$, its presence does have an impact on analytic structure of its sum (or $\frac{1}{\mu(1-\alpha)}$-sum, if $F$ diverges).

\begin{example}

Let  $f(t)$ be the entire function $\op{e}^{t}$. Assuming $\mu=1$ and using this $f$ and the ground state $\psi_{\alpha,\, \beta,\,\gamma}$, \cref{eq:rigorous_umbral_image_psi_alpha_beta} becomes:
\begin{dmath}\label{eq:ML_3par__series}
	F(\zeta) = {
		\op{e}^{\zeta\uo} \,[\psi_{\alpha,\, \beta,\,\gamma}] =
		\sum_{r = 0}^\infty \frac{(\gamma)_r}{r!\,\Gamma(\alpha r + \beta)}\, \zeta^r \condition*{\alpha >0,\; \beta, \gamma \in \bbb{C}}
	}
\end{dmath},
which coincides with the $(\frac{1}{\alpha},\,\beta)$-Borel transform of 
\begin{dmath*}
		\tilde{h}(t) = {
			t\,h(t) = 
			t\sum_{r = 0}^\infty \frac{(\gamma)_r}{r!}\,t^r
		}
\end{dmath*},
converging to $t/(1-t)^\gamma$ in an open disk centred in the origin.

Since $h \in \bbb{C}\llb t\rrb_0$, \cref{eq:gevrey_nested} and \cref{eq:k-gevrey-theorem} ensure that $F$ converges for any $\zeta \in \bbb{C}$ and any $\alpha > 0$. \Cref{eq:ML_series} is in fact the Maclaurin expansion of the three-parameters Mittag-Leffler function $E_{\alpha,\,\beta,\,\gamma}(\zeta)$, an entire function of $\zeta$ also known as Prabhakar function \cite{tP71}.

The $(\frac{1}{\alpha},\,\beta)$-Laplace transform of this function yields as expected:
\begin{dmath*}
	\mathcal{L}_{\frac{1}{\alpha},\,\beta}[E_{\alpha,\,\beta,\,\gamma}](t) = \frac{t}{(1-t)^\gamma} 
\end{dmath*}.

In conclusion, the following umbral identity holds for any $\zeta \in \bbb{C}$:
\begin{dmath*}
	E_{\alpha,\,\beta,\,\gamma}(\zeta) = \op{e}^{\zeta\uo}[\psi_{\alpha,\,\beta,\,\gamma}]
\end{dmath*}.

\end{example}
\vskip 0.2cm

\section{A case study: Gaussian trigonometric functions}\label{sec:gaussian_trig}

As an interesting application of the umbral framework outlined in the previous sections, we investigate a new umbral formulation of the so-called Gaussian trigonometric functions, namely:
\begin{dgroup}\label[pluralequation]{eq:gaussian_trig}
	\begin{dmath}\label{eq:gaussian_exp}
		\op{exp}_G(\zeta) \equiv {\op{w}(\zeta) \defeq \op{e}^{-\zeta^2}(1+i\,\op{erfi}(\zeta))}
	\end{dmath},
	\begin{dmath}\label{eq:gaussian_cos}
		\cos_G(\zeta) \defeq \op{e}^{-\zeta^2}
	\end{dmath},
	\begin{dmath}\label{eq:gaussian_sin}
		\sin_G(\zeta) \defeq \op{e}^{-\zeta^2} \op{erfi}(\zeta)
	\end{dmath}.
\end{dgroup}

The Gaussian exponential (a.k.a. Faddeeva) function \cite{vnF54}, is a fundamental special function in mathematical physics, applied mathematics, and computational science. It plays a central role in the evaluation of convolution integrals, spectral line shapes, resonance broadening and kinetic plasma theory, where it is usually referred to as plasma dispersion function. The spectrum of its applications is very wide and pervasive across disciplines.

In \cite{gD23}, the following umbral formulation of these functions was provided in terms of the umbral state function $\phi$ (see \cref{eq:umbral_state_phi}):
\begin{dgroup}\label[pluralequation]{eq:umbral_trig_phi}
	\begin{dmath}\label{eq:umbral_exp}
		\op{exp}_G(\zeta) = \frac{1}{1 - i \zeta \uo^{1/2}}\,[\phi]
	\end{dmath}
	\begin{dmath}\label{eq:umbral_cos_phi}
		\cos_G(\zeta) = \frac{1}{1+\zeta^2\uo}\,[\phi]
	\end{dmath},
	\begin{dmath}\label{eq:umbral_sin_phi}
		\sin_G(\zeta) = \frac{\zeta \uo^{1/2}}{1+\zeta^2\uo}\,[\phi]
	\end{dmath},
\end{dgroup}

On the other hand, we note that, by exploiting the following identity:
\begin{dmath*}
	\op{exp}_G(\zeta) \equiv {\op{w}(\zeta)= 
		E_{\frac{1}{2},\,1}(i\zeta)
	}
\end{dmath*},
it is possible to define alternative umbral representations for the Gaussian exponential, sine and cosine functions that emphasise their trigonometric-like nature. Namely, we can establish the simple umbral identities (cfr. \cref{eq:faddeeva_function}:
\begin{dgroup}\label[pluralequation]{eq:umbral_trig}
	\begin{dmath}\label{eq:umbral_exp}
		\op{exp}_G(\zeta) = \op{e}^{i\zeta\uo}\,[\lambda]
	\end{dmath},
	\begin{dmath}\label{eq:umbral_cos}
		\cos_G(\zeta) = \cos(\zeta\uo)\,[\lambda]
	\end{dmath},
	\begin{dmath}\label{eq:umbral_sin}
		\sin_G(\zeta) = \sin(\zeta\uo)\,[\lambda]
	\end{dmath},
\end{dgroup}
where $\lambda$ is the already introduced ground state function
\begin{dmath}\label{eq:lambda}
	\lambda \defeq {
		\psi_{\frac{1}{2},\, 1} =
		\frac{\Gamma(1+t)}{\Gamma(1 +t/2)}
	}
\end{dmath}.

So, the ordinary exponential, sine and cosine functions are respectively the $\lambda$-umbral images of the Gaussian exponential, sine and cosine.

\Cref{eq:umbral_trig} can be fruitfully used in many different applicative contexts. E.g., by using \cref{eq:umbral_cos} we can deduce a simple $\lambda$-umbral image for the error function:
\begin{dmath*}
	\op{erf}(z) = {
		\frac{2}{\sqrt{\pi}}\int_0^z \mathrm{d}\zeta \cos_G(\zeta) =
		\frac{2}{\sqrt{\pi}}\int_0^z \mathrm{d}\zeta \cos(\zeta\uo)\,[\lambda] =
		\frac{2}{\sqrt{\pi}}\uo^{-1}\sin(z\uo)\,[\lambda]
	} =
	\frac{2}{\sqrt{\pi}}\sum_{r=0}^\infty\,\frac{(-1)^rz^{2r+1}}{(2n+1)n!}
\end{dmath*}.
So:
\begin{dmath*}
	\op{erf}(\zeta) = \frac{2}{\sqrt{\pi}}\frac{\sin(\zeta\uo)}{\uo}\,[\lambda]
\end{dmath*}.

As another example, the formula
\begin{dmath*}
	\op{e}^{i x u} = -\frac{i}{\pi}\mathcal{P}\int_\infty^\infty \mathrm{d}\xi\,\frac{\op{e}^{i\xi u}}{\xi - x} \condition*{u>0}
\end{dmath*},
where the symbol $\mathcal{P}$ indicates the Cauchy principal value of the integral, directly translates into the umbral identity
\begin{dmath*}
	\op{exp}_G(\zeta) = {
		\op{e}^{i x \uo} \,[\lambda] = 
		-\frac{i}{\pi}\mathcal{P}\int_\infty^\infty \mathrm{d}\xi\,\frac{\op{e}^{i\xi \uo}}{\xi - x}\,[\lambda] = 
		-\frac{i}{\pi}\mathcal{P}\int_\infty^\infty \mathrm{d}\xi\,\frac{\exp_G(\xi)}{\xi - x}
	}
\end{dmath*}.
Using the functional identity $\op{exp}_G(\zeta) = \cos_G(\zeta) + i\sin_G(\zeta)$, 
one is eventually led to  Kramers-Kronig relations among the Gaussian sine and cosine:
\begin{dgroup}\label[pluralequation]{eq:umbral_trig_KK}
	\begin{dmath}\label{eq:umbral_cos_KK}
		\cos_G(z) = {
			\frac{1}{\pi}\mathcal{P}\int_\infty^\infty \mathrm{d}\xi\,\frac{\sin(\xi \uo)}{\xi - x}\,[\lambda] =
			\frac{1}{\pi}\mathcal{P}\int_\infty^\infty \mathrm{d}\xi\,\frac{\sin_G(\xi)}{\xi - x}
		}
	\end{dmath},
	\begin{dmath}\label{eq:umbral_sin_KK}
		\sin_G(z) = {
			-\frac{1}{\pi}\mathcal{P}\int_\infty^\infty \mathrm{d}\xi\,\frac{\cos(\xi \uo)}{\xi - x}\,[\lambda] =
			-\frac{1}{\pi}\mathcal{P}\int_\infty^\infty \mathrm{d}\xi\,\frac{\cos_G(\xi)}{\xi - x}
		}
	\end{dmath}.
\end{dgroup}

The umbral versions of $n$-th order derivatives of the Gaussian trigonometric functions are obtained by directly deriving \cref{eq:umbral_trig}:
\begin{dgroup}\label[pluralequation]{eq:umbral_trig_derivatives}
	\begin{dmath}\label{eq:umbral_exp_n}
		\exp_g^{(n)}(x) = (i\uo)^n\op{e}^{ix\uo}\,[\lambda]
	\end{dmath},
	\begin{dmath}\label{eq:umbral_cos_n}
		\cos_G^{(n)}(x) = \uo^n \cos\left(x\uo + \frac{\pi n}{2}\right)\,[\lambda] =
			\cos\left(\frac{\pi n}{2}\right)\uo^n\cos(x\uo)\,[\lambda]  - \sin\left(\frac{\pi n}{2}\right) \uo^n \sin(x\uo)\,[\lambda]
	\end{dmath},
	\begin{dmath}\label{eq:umbral_sin_n}
		\sin_G^{(n)}(x) = \uo^n \sin\left(x\uo + \frac{\pi n}{2}\right)\,[\lambda] =
			\cos\left(\frac{\pi n}{2}\right)\uo^n\sin(x\uo)\,[\lambda]  + \sin\left(\frac{\pi n}{2}\right) \uo^n \cos(x\uo)\,[\lambda] 
	\end{dmath}.
\end{dgroup}
From \cref{eq:umbral_exp_n} we immediately derive the formula linking the Gaussian exponential a.k.a. Faddeeva function to its first derivative \cite[\href{https://dlmf.nist.gov/7.10.E2}{(7.10.2)}]{NIST:DLMF}, namely:
\begin{dmath*}
	\op{w}'(x)	 = {
		i\uo\op{e}^{ix\uo}\,[\lambda] =
		\sum_{r=0}^\infty \frac{(i)^{r+1}x^r}{r!}\frac{\Gamma(r+2)}{\Gamma\left(\frac{r}{2}+\frac{3}{2}\right)} =
		\frac{2i}{\sqrt{\pi}} - 2x\op{w}(x)
	}
\end{dmath*}.

We can use \cref{eq:umbral_cos_n} to recast the defining relation of the Hermite polynomials,
\begin{dmath*}
	H_n(x) \op{e}^{-x^2}	 = (-1)^n \partial_x^n \op{e}^{-x^2}
\end{dmath*},
into the simple umbral form:
\begin{dmath*}
	H_n(x) \op{e}^{-x^2}	 = {(-1)^n \cos_G^{(n)}(x) = (-\uo)^n \cos\left(x\uo + \frac{\pi n}{2}\right)\,[\lambda] }
\end{dmath*}.

By exploiting the classical Dirichlet integral
\begin{dmath}
	\int_{-\infty}^\infty \mathrm{d}x\,\frac{\sin(u x)}{x} = \pi \sgn(u) \condition*{u \in \bbb{R}}	
\end{dmath},
we immediately obtain a simple result for the non trivial integral\footnote{We remark that $c\,[\varphi] = c \uo^0 \,[\varphi] = c \varphi(0)$, for any constant $c$ and ground state $\varphi$. In the present case, $\lambda(0) = 1$.}:
\begin{dmath}
	\int_{-\infty}^\infty \mathrm{d}x\,\frac{\op{e}^{-x^2} \op{erfi}(x)}{x} = {
	\int_{-\infty}^\infty \mathrm{d}x\,\frac{\sin_G(u x)}{x} =
		\int_{-\infty}^\infty \mathrm{d}x\,\frac{\sin(x \uo)}{x}\,[\lambda] =
		\pi \,[\lambda] = \pi
	}
\end{dmath}.

\subsection{Connection with the Fourier transform}

The specific form of the $\lambda$-umbral image of the Gaussian exponential, see \cref{eq:umbral_exp}, suggests an intriguing  connection with the theory of the Fourier transform. Indeed, we have that formally:
\begin{dmath}\label{eq:umbral_Fourier_transform}
	\int_{-\infty}^\infty \mathrm{d}x\,f(x) \exp_G(-kx) = {
		\int_{-\infty}^\infty \mathrm{d}x\,f(x) \op{e}^{-i x k\uo}\,[\lambda] =
		\mathcal{F}[f](k\uo)\,[\lambda]
	}
\end{dmath},
where $\mathcal{F}[f]$ denotes the Fourier transform\footnote{
	We adopt the following conventions for the Fourier transform and its inverse:
	\begin{dmath*}
		\hat{f}(k) = {\mathcal{F}[f](k) = \int_{-\infty}^\infty \mathrm{d}x\,f(x) \op{e}^{-ikx}, \quad\quad f(x) = \mathcal{F}^{-1}[\hat{f}](x) = \int_{-\infty}^\infty \frac{\mathrm{d}k}{2\pi}\,f(x) \op{e}^{ikx}} 
	\end{dmath*}.
} 
of $f$.

\noindent If we interpret the first integral as the "Gaussian Fourier transform" of $f$ and introduce the suggestive notation:
\begin{dmath}\label{eq:gaussian_Fourier_transform}
	\mathcal{F}_G[f] = \int_{-\infty}^\infty \mathrm{d}x\,f(x) \exp_G(-kx) 
\end{dmath},
\cref{eq:umbral_Fourier_transform} can be rewritten as
\begin{dmath}\label{eq:_umbral_gaussian_Fourier_transform}
	\mathcal{F}_G[f](k) = \mathcal{F}[f](k\uo)\,[\lambda] 
\end{dmath}.
 In other words, the Gaussian Fourier transform of $f$ is the formal series whose $\lambda$-umbral image is the ordinary Fourier transform of $f$.

Since the convergence conditions of $f$ and $\mathcal{F}[f]$ can differ, \cref{eq:_umbral_gaussian_Fourier_transform} can yield either a convergent formal series -- i.e. ultimately a function -- or a divergent but resummable formal series -- i.e. the asymptotic expansion of a function in an open sector of the complex plain. With this caveat, \cref{eq:umbral_Fourier_transform} provides a simplified method to perform complex calculations -- e.g. integrals involving the Gaussian exponential -- by exploiting known Fourier transforms. For example, the integral:
\begin{dmath*}
	\int_{-\infty}^\infty \mathrm{d}x\op{e}^{-ax^2} \exp_G(-x) = {
		\int_{-\infty}^\infty \mathrm{d}x\op{e}^{-(a+1)x^2}(1 + \op{erfi}(-x)) = \sqrt{\frac{\pi}{a+1}} \condition*{\op{Re}a>0}
	}
\end{dmath*}
can be calculated as the value at $k=1$ of the Gaussian Fourier transform of $\op{e}^{-ax^2}$, namely:
\begin{dmath}\label{eq:ex_umbral_gaussian}
	\int_{-\infty}^\infty \mathrm{d}x\op{e}^{-ax^2} \op{e}^{-i x \uo}\,[\lambda] = {
		\sqrt{\frac{\pi}{a}} \op{e}^{-\frac{\uo^2}{4a}}\,[\lambda] = 
		\sqrt{\frac{\pi}{a}}\sum_{r=0}^\infty \frac{(-1)^r\,(2r)!}{(4a)^r (r!)^2} =
		\sqrt{\frac{\pi}{a+1}}
	}
\end{dmath}.

We provide in the following some other examples of application of the method.

\begin{example}

Let us consider the integral
\begin{dmath*}
	I = {\int_{-\infty}^\infty \mathrm{d}x\op{e}^{-|x|} \exp_G(-x) = \sqrt{\pi}\op{e}^{\frac{1}{4}}(1 - \op{erf}(1/2))}
\end{dmath*}.
We note that $I$ corresponds to $\mathcal{F}_G[\op{e}^{-|x|}](k)$ evaluated in $k=1$.
The method just outlined yields in this case:
\begin{dmath*}
	F(k) \equiv {
	\mathcal{F}_G[\op{e}^{-|x|}](k) = 
		\mathcal{F}[\op{e}^{-|x|}](k\uo)\,[\lambda] =
		\frac{2}{1+k\uo^2}\,[\lambda]
	}
\end{dmath*}.
We have already encountered this umbral identity in \cref{ex:7}, where we established that the resulting formal series is divergent 1-Gevrey and 1-summable to the function:
\begin{dmath*}
	\op{S}[F](k) = 
		\sqrt{\frac{\pi}{k}}\op{e}^{1/4k} \op{erfc}(1/2k^{1/2})\condition*{k \in V=S(\Theta),\;\Theta=2\pi}
\end{dmath*}.
We eventually obtain:
\begin{dmath*}
	I = {\op{S}[F](1) = \sqrt{\pi}\op{e}^{\frac{1}{4}}(1 - \op{erf}}(1/2))
\end{dmath*},
which coincides with the result obtained by directly solving the integral\footnote{Note that $\op{erfc}(\zeta) \defeq 1 - \op{erf}(\zeta)$.}.

\end{example}
\vskip 0.2cm

\begin{example}

\begin{dmath*}
	\mathcal{P}\int_{-\infty}^\infty \frac{\exp_G(-x)}{x} = {
		\mathcal{F}[1/x](\uo)\,[\lambda] =
		-i\pi\,[\lambda] =
		-i\pi
	}
\end{dmath*}.

\end{example}
\vskip 0.2cm

\begin{example}

Let us consider the integral
\begin{dmath*}
	I_n = {
		\int_{-\infty}^\infty \mathrm{d}x\op{e}^{-x^2/2} H_n(x) \exp_G(-x) =
		\mathcal{F}_G[\op{e}^{-x^2/2} H_n(x)](1)
	}
\end{dmath*},
where $H_n(x)$ is the Hermite polynomial of order $n \in \bbb{N}$.
The Gaussian Fourier transform yields the neat result
\begin{dmath}\label{eq:umbral_hermite_gf}
	\mathcal{F}_G[\op{e}^{-x^2/2} H_n(x)](k) = {
		(-i)^{\,n}\sqrt{2\pi}\op{e}^{-k^2 \uo^2/2} \, H_n(k\uo)\,[\lambda]
	}
\end{dmath}.
In order to verify the convergence of \cref{eq:umbral_hermite_gf}, we note that, for each $n$, it can be written as a linear combination, with given coefficients, of terms of the form:
\begin{dmath*}
	k^m\op{e}^{-k^2\uo^2/2} \, \uo^m\,[\lambda] \condition*{m \leq n}
\end{dmath*},
hence it suffices to analyse this particular case. We have:
\begin{dmath*}
	k^m\op{e}^{-k^2\uo^2/2} \, \uo^m\,[\lambda] = {
		k^m\sum_{r=0}^\infty\,\frac{(-1)^r}{2^r r!}\frac{\Gamma(2r + m + 1)}{\Gamma(r + 1 + m/2)}\,k^{2r}
	}
\end{dmath*}.
This series converges to the entire function of $k$:
\begin{dmath*}
	a_m(k) = \frac{k^m (2 k^2+1)^{-\frac{m-1}{2}} \,\Gamma (m+1)}{\Gamma (1+m/2)}
\end{dmath*}, whence:
\begin{dmath}\label{eq:coefficients}
	a_m(1) = \frac{3^{-\frac{m-1}{2}}\, \Gamma (m+1)}{\Gamma (1+m/2)}
\end{dmath}.
In conclusion, \cref{eq:umbral_hermite_gf} converges to an entire function of $k$ for any $n \in \bbb{N}$. When evaluated in $k=1$, this function returns the desired value of $I_n$.

\noindent Direct calculations of e.g. $I_1$ and $I_3$, carried out by using \cref{eq:coefficients}, give:
\begin{dmath*}
	I_1 = {
		\int_{-\infty}^\infty \mathrm{d}x\op{e}^{-x^2/2} (2x) \exp_G(-x) =
		-\frac{4i\sqrt{2}}{3}
	} 
\end{dmath*},
\begin{dmath*}
	I_3 = {
		\int_{-\infty}^\infty \mathrm{d}x\op{e}^{-x^2/2} (8x^3 -12x) \exp_G(-x) =
		-\frac{8 i \sqrt{2}}{9}
	} 
\end{dmath*},
which correspond to the results obtained by directly solving the integrals with the aid of Mathematica\texttrademark.

\end{example}
\vskip 0.2cm

\section{Conclusions}\label{sec:conclusions}

We have demonstrated that the indicial umbral formalism, developed as a computational tool in the framework of special function theory, can be profitably transposed in the wider context of formal power series $\bbb{C}\llb t \rrb$, by means of a suitable redefinition of the umbral operator $\uo$ as a functional in the ground-state subalgebra of analytically convergent formal series $\bbb{C}\{t\} \subset \bbb{C}\llb t \rrb$.

For a carefully defined subset of ground states $\varphi$, wide enough to comprise most of the interesting cases investigated so far, this transposition enables to make rigorous sense of umbral identities of the form $F(\zeta) = f(\zeta \uo^\mu)[\varphi]$, reinterpreted as well-defined formal series $F \in \bbb{C}\llb \zeta \rrb$ for any $f \in \bbb{C}\llb t \rrb$.

In order to make contact with concrete applications of the formalism, which require $F \in \bbb{C}\{\zeta\}$, we have imposed the additional condition $f \in \bbb{C}\{t\}$.

Even with this restriction on $f$, the formal series $F(\zeta) = f(\zeta \uo^\mu)[\varphi]$ can still be divergent for certain values of the parameters entering the ground state $\varphi$. We have thus analysed its convergence conditions, by applying well-established results of the Gevrey classification of formal series.

This has eventually enabled to determine some general criteria for the convergence of $F$. In particular, we have established that, even when $F$ is divergent of Gevrey level $k>0$, it can possibly be resummed by means of the generalised $k$-Borel-Laplace summation procedure. These results have been verified in a series of explicit examples and applied to the analysis of a new promising umbral formulation of the Gaussian trigonometric functions.

\begin{table}[h!]
\centering
\renewcommand{\arraystretch}{1.1}
\small
\begin{tabular}{|p{3cm}|p{4cm}|p{5cm}|}
\hline
\textbf{} & \textbf{Classical umbral calculus} & \textbf{The present approach} \\
\hline
\textbf{Ground state $\varphi$} & Arbitrary formal series; mainly combinatorial & Structured functions (e.g., involving $\Gamma$, Pochhammer); tied to summability \\
\hline
\textbf{Analytic control} & Minimal; purely symbolic & Detailed asymptotics and convergence via Gevrey theory \\
\hline
\textbf{Goal} & Symbolic generation of identities and sequences & Rigorous construction of functions, possibly via resummation techniques \\
\hline
\textbf{Borel-Laplace theory} & Not used & Central tool for reconstructing functions from divergent series \\
\hline
\textbf{Operator $\uo^{\,\mu}$} & Symbolic shift or lowering operator & Realized as a functional; carefully analysed convergence \\
\hline
\textbf{Gevrey classes} & Not present & Essential to classify and resum formal series \\
\hline
\textbf{Function-analytic viewpoint} & Absent & Strong: involves growth, singularities, domains of analyticity \\
\hline
\end{tabular}
\caption{\small Key Differences with Roman's Umbral Calculus}
\end{table}

\vskip 0.3cm
We conclude with a brief critical comparison of the umbral framework proposed in this work with Roman's classical umbral calculus. The results of the comparison are summarised in table 1. 

Although both approaches share many points of contact, they are essentially different in scope and aimed at different objectives. In Roman's theory expressions of the form $f(\uo)[\varphi]$ are interpreted algebraically. These expressions serve as combinatorial tools to represent Sheffer sequences, linear functionals, and other symbolic constructs.

By contrast, in the present approach  umbral expressions such as $f(\zeta \uo^\mu)[\varphi]$ are reinterpreted 
within the rigorous framework of formal power series, and their convergence analysed in terms of Gevrey classes and Borel-Laplace summation. Here, the ground state $\varphi \in \bbb{C}\{t\}$ is typically the Maclaurin expansion of a structured function (e.g., a ratio of Gamma functions or Pochhammer symbols), and the operator $\uo^\mu$ is eventually understood as a functional on the subalgebra of convergent formal series, constrained to produce a convergent result. Stated differently, our formalism maintains symbolic generality, but is anchored in analytic properties. In particular:
\begin{itemize}
  \item formal series are assigned to Gevrey classes based on the growth of their coefficients;
  \item when divergent, these series are interpreted through Borel-Laplace resummation, yielding analytic functions on sectors of the complex plane;
  \item the umbral expression $f(\zeta \uo^\mu)[\varphi]$ is thus not merely symbolic: its convergence and analytic structure depend on the interplay between the function $f$ and the ground state $\varphi$, which determine the Gevrey type of the associated formal series.
\end{itemize}
So, while Roman's umbral calculus is algebraic and combinatorial in spirit, the present "analytic" umbral formalism endows the underlying algebraic structure with tools from asymptotic analysis and summability theory, yielding a powerful framework for constructing and interpreting special functions from divergent formal data.

\vskip 0.3 cm
In conclusion, the proposed framework provides a solid foundation for a wide subset of the indicial umbral methods developed so far, by also enabling to make sense of some apparently diverging results. 

\noindent Work is in progress to extend the framework by including further classes of admissible ground states, as well as different types of resummation procedures (such as Borel-Le Roy resummation), and apply these ideas to the investigation of some fundamental problems of modern theoretical physics, such as certain aspects of renormalisation theory \cite{gP79}.

\vskip 0.3 cm

\noindent The author is grateful to Dr. G. Dattoli for reading the manuscript and providing many illuminating suggestions. 

\bibliographystyle{plain}
\bibliography{bibliography}


\end{document}